\documentclass[10pt]{article}
\usepackage{amssymb}
\usepackage[utf8]{inputenc}
\usepackage[english]{babel}
\usepackage{graphicx}
\usepackage{fancyhdr}
\usepackage{multicol}
\usepackage{multirow}
\usepackage{amsmath, amsthm, amssymb, amsfonts}
\usepackage{float}
\usepackage{enumerate}
\usepackage{cancel}
\usepackage{tikz}
\numberwithin{equation}{section}
\usepackage[hidelinks]{hyperref}
\hypersetup{colorlinks=true,linkcolor=blue}
\newcommand{\R}{\mathbb R}

\def\R{\mathbb{R}}

\def\> {\rightarrow}

\def\>{\rightarrow}
\def\S1{\mathbb S^1}

\newtheorem{Definicion}{Definici\'on}[section]
\newtheorem{example}[Definicion]{Example}

\newtheorem{teo}[Definicion]{Theorem}

\newtheorem{propo}[Definicion]{Proposition}
\newtheorem{lem}[Definicion]{Lemma}
\newtheorem{coro}[Definicion]{Corollary}

\newcommand{\fin}{\hfill$\blacksquare$}

\usepackage{fancyhdr}
\title{Periodic solutions for one-dimensional nonlinear nonlocal problem with drift including singular nonlinearities}
\author{Lisbeth Carrero and Alexander Quaas}
\date{}
\begin{document}
\maketitle

\abstract{In this paper, we prove existence results of a one-dimensional periodic solution to equations with the fractional Laplacian of order $s\in(1/2,1)$, singular nonlinearity, and gradient term under various situations, including nonlocal contra-part of classical Lienard vector equations, as well other nonlocal versions of classical results know only in the context of second-order ODE.  Our proofs are based on degree theory and Perron's method, so before that we need to establish a variety of priori estimates under different assumptions on the non-linearities appearing in the equations. Besides, we obtain also multiplicity results in a regime where a priori bounds are lost and bifurcation from infinity occurs.}

\section{Introduction}
Second-order ordinary differential equations (ODE) and systems are since Newton's second law of motion, one of the most study equations in mathematics and physics under many different situations. They also play a crucial role in the study of linear and nonlinear PDE's. In this work, we want to extend some results for second-order ODE found in the work of Mawhin, see \cite{mawhin} to equations with the one-dimensional fractional laplacian instead of laplacian. In particular, we will prove the existence of periodic solutions to some basic classical model like Lienard \cite{Lienard}(1928), Forbat, F., Huaux \cite{Forbat}(1962) Lazer-Solimini \cite{Lazer} (1987)
for equations involving the fractional laplacian, see also other references in \cite{mawhin}.\\

Nonlocal operator with singular kernels, in particular, the fractional Laplacian has received much attention in recent years; this motived by many apply models in biology, physics, chemistry, finance, marine foreign, etc., where the underlining phenomena are governed by anomalous diffusion, connected with Levy flights where the fractional Laplacian appears naturally. As an example of applied phenomena where this type of operator appears we mention \cite{Hu}, \cite{Metzler}, \cite{V1}, \cite{V2}, \cite{fisica} and a more mathematical review on the topic can be found in \cite{val}, see also the big list of reference in all these works.
Notice that periodic patterns are naturally expected in many of these applied phenomena.\\

The fractional Laplacian operator can be defined for $s\in (0,1)$ via its multiplier $|\xi|^{2s}$ in Fourier space, notice $s=1$ corresponds to the Laplacian. It can also be defined by the formula
\begin{equation}\label{fracLapl}
(-\Delta)^s u(x) = {C_{N, s} \mathrm{P.V.} \int_{\mathbb{R}^N} \frac{u(x) - u(z)}{|x - z|^{N + 2s}}dz},
\end{equation}
here $C_{N, s} > 0$ is a well-known normalizing constant and P.V stands for the principal value, see for example \cite{val}.

From the mathematical point of view, the equations involving the fractional Laplacian that we will study here required some new different technics since many tools from ODE such as energy method, integral factors, and other key elements do not remain valid in the nonlocal case, and therefore we need  different method and argument for the non-local case.\\

Before describing our main results, notice that the type of equations studied in \cite{mawhin} are of the form
 \begin{equation}\label{m1}
u''(t) + f(u(t))u'(t) + g(u(t))= e(t) \;\;\;\; t\in\mathbb{R},
\end{equation}
and
 \begin{equation}\label{m2}
u''(t) + cu'(t) - g(u(t))= e(t) \;\;\;\; t\in\mathbb{R}.
\end{equation}
Mawhin was mainly interested in singular nonlinearities, that is, he supposes that $g$ becomes unbounded near the origin; the $+$ sign indicates that the particles have opposite charges, while the $-$ sign indicates that the particles have the same charge.
 So, it is said that the equation with the $+$ sign (resp. $-$) has an attractive singularity (repulsive resp.). We simply talk about the attractive and repulsive case.\\

Thus, want studied equation of the form
 \begin{equation}\label{intro}
(\triangle)^{s}u(t) + f(u(t))u'(t) + g(u(t))= e(t) \;\;\;\; t\in\mathbb{R},
\end{equation}
and
\begin{equation}\label{ecuag2}
  (\triangle) ^{s}u(t) + cu'(t) - g(u(t)) = e(t)  \;\;\;\; t\in\mathbb{R},
\end{equation}
where $(\triangle)^{s}u(t):= - (-\triangle)^{s}u(t)$ and $s\in (1/2,1)$ from now on.\\

The main aim of the present paper, is to establish existence of periodic solutions to equations $(\ref{intro})$ and $(\ref{ecuag2})$.

In the study of (\ref{intro}) we suppose that $f:(0, + \infty) \to \mathbb{R}$, $g:(0, + \infty) \to \mathbb{R}$ are $C^{\alpha}((0, + \infty))$ and $e\in C^{\alpha}(\mathbb{R})$ with $\alpha \in (0,1)$. Without loss of generality, we assume that we are searching for $2\pi$-periodic solutions to simplify the discussion, therefore, we will assume that $e$ is also $2\pi$-periodic.
Moreover,  we will use denotation
\begin{equation*}
  \bar e:=\dfrac{1}{2\pi}\int_0^{2\pi}e(t)dt,
\end{equation*}
its mean value.\\
Observe that the term with $u'$ is sometimes called drift term and corresponds to some "transport" or "friction" in some of the models.

Before giving our results, let us mention that periodic solutions are studied in \cite{quaas}, \cite{Wei} without drift term and regular nonlinearity, see also \cite{du}, \cite{gui}, \cite{juan}.

Other types of periodic problems related to the spectral fractional Laplacian can be found in \cite{ambrosio, ambrosio1, ambrosio3, ambrosio4, ambrosio2}.

Now, we will give our first existence result in the case of attractive singularity.

\begin{teo}\label{teo}
Assume that the function $g:(0, + \infty) \to \mathbb{R}$ is such that the following conditions hold:
\begin{enumerate}
 \item $g(t) \to + \infty$ as $t\to 0^{+}$.
 \item $\limsup\limits_{t \to + \infty} g(t)<\bar{e}$.
\end{enumerate}
Then equation (\ref{intro}) has at least one $2\pi$-periodic positive classical solution.
\end{teo}

This theorem will be proved by using Perron's method. One of the main difficulty here is to find a periodic super-solution, this is obtained by solving a semi-linear problem of Lienard type. In fact, we have the following results for nonlocal Lienard type equations.

\begin{propo}\label{l1}
  Let $f\in C^{\alpha}(0, + \infty)$ and $w\in C^{\alpha}_{2\pi}(\mathbb{R})$, then
\begin{equation}\label{le}
(\triangle)^{s}u(t) + f(u(t))u'(t)=w(t),
\end{equation}
has at least one $2\pi$-periodic classical solution $u$ if, only if $\bar{w}=0$. Moreover,
\begin{equation*}
  \|u \|_{C^{2s + \alpha}(\mathbb{R})} \leq M\|w\|_{C^{\alpha}(\mathbb{R})},
\end{equation*}
for some positive constant $M$.
\end{propo}

To prove the necessary condition, we will use a basic property of  fractional Laplacian for periodic functions, see Lemma \ref{perizero}. To prove the existence
will use the Schauder type estimates which are obtained by $H^s$ estimate that gives by the assumption $s\in(1/2, 1)$ a $C^\alpha$ estimate then with the help of interpolation inequality in H\"older space we are able to manage the drift term, so we get our Schauder type estimates.\\

This type of results can be generalized to the nonlocal Lienard vector equation inspired by \cite{principal}.
More precisely, we will study vector equations of the form:

\begin{equation}\label{sistema}
 L(u(t)) + \frac{d}{dt}( \nabla H(u(t))) + Au(t)= e(t),
\end{equation}
where
\begin{equation*}
L(u(t))= \left( \begin{array}{c}
 (\triangle) ^{s_{1}}u_{1}(t) \\
 \vdots \\
 (\triangle)^{s_{n}} u_{n}(t)
 \end{array}\right) \;\; s_{i}\in\left(\frac{1}{2},1\right)\;\; i= 1, \cdots, n,
\end{equation*}
with $H: \mathbb{R}^{n} \to \mathbb{R}$ is $C^{2, \alpha}( \mathbb{R}^{n})$, $A$ is a $n\times n$- matrix and $e\in C^{\alpha}(\mathbb{R}, \mathbb{R}^{n})$ and $e$ is also $2\pi$-periodic.\\
Looking for existence results, we give a sufficient condition for (\ref{sistema}) to have a $2\pi$-periodic solution.

\begin{teo}\label{teosistem}
If $\bar{e}\in Im A$ and
\begin{equation*}
 M= \sup\limits_{|y| = 1} \big< Ay, y\big> \;\; < 1,
\end{equation*}
then (\ref{sistema}) has at least one $2\pi$-periodic classical solution.
\end{teo}

Now we describe our second main results concerning equation $(\ref{ecuag2})$ that is in the case of repulsive singular nonlinearity.

For that, we suppose that $c>0$, $e\in C^{\alpha}(\mathbb{R})$, $e$ is $2\pi$-periodic with $\bar{e} > 0$ and $g \in C^{\alpha}(0, + \infty)$ is a given function satisfying the following conditions:

\begin{enumerate}[(G1)]
 \item $\limsup\limits_{t \to + \infty} [g(t) + \bar{e}]<0$,
 \item $\lim\limits_{t \to  0^{+}} g(t) = + \infty$ and $g$ is monotone near zero,
 \item There exists $\epsilon>0$ such that $g(\tau)^{2s-2-\epsilon}\int\limits_{\tau}^{1} g(t) dt \to + \infty$ as $\tau \to 0^{+}$,
 \item $g(t)\geq -at -b$ for some $a>0$ and $b\geq0$ and all $t>0$.
\end{enumerate}

Now we are in a position to give our second main theorem.

\begin{teo}\label{teoulti}
Assume that the conditions $(G1)$, $(G2)$,$(G3)$ and $(G4)$ are fulfilled, then equation (\ref{ecuag2}) has at least one $2\pi$-periodic positive classical solution.
\end{teo}

Here, is important to mention the main difficulties with respect to the local case ($s=1$) where a priori uniform bounds hold for $u$, $u'$. But in the nonlocal case we just found  a $L^{2}$ bound for $u'$.

Besides that, the main idea to get  the a priori bound is to prove a new general energy type identity for periodic solution (see Lemma  \ref{lemint}) that gives the formula for $$\int_a^b (\Delta^s)u(t) u'(t)dt=I(u,a,b),\quad \mbox{for all }\quad a,b \in \R,$$ inspired in \cite{simetria} that give related formula in the half-line.
Here $I(u,a,b)$ is the right side given in Lemma  \ref{lemint}. This identity together with a local regularity theory (Harnack inequality, and Schauder type estimates) at the maximum point of a solution will help us to estimate $I(u,a,b)$, then using $(G3)$ we are able to find the a priori lower bound for our periodic solutions, for more details see section $4$. The rest of the proof is base on degree type arguments.\\

As a by product of Theorem \ref{teoulti} and bifurcation from infinity, this again inspired from the works \cite{mawhin} and \cite{schmitt}, we find multiplicity results in a regime where a priori bounds are lost, see Theorem \ref{bifu}, below.\\

Let finish this introduction by mentioning that some reference of previous results for fractional equation with singular nonlinearities in a different situation from the studied here can be found in \cite{barrios, canino, gui}.

The paper is organized as follows. In Section $2$, is devoted to proving Proposition $\ref{l1}$ which is the semi-linear case of $(\ref{intro})$ and its generalization Theorem \ref{teosistem}. Section $3$ and $4$ are devoted to establishing some results for the existence of at least one $2\pi$- periodic positive solution for (\ref{intro}) and (\ref{ecuag2}) respectively. Finally, the solution multiplicity results will be presented in Section $5$.

\section{The semi-linear periodic problems}

In this section, first we will prove some preliminary results that will facilitate the proof of Theorem \ref{teosistem} and Proposition $\ref{l1}$.
Notice that by Proposition $\ref{l1}$  of \cite{quaas} the fractional Laplacian of a smooth, $2\pi$-periodic function reduces to
\begin{equation}\label{lui}
    (\triangle)^{s}u(t) =-\mathcal{L}u(t) \;\;\; t\in(0,2\pi),
\end{equation}
where
\begin{equation}\label{L}
   \mathcal{L}u(t)=  \int\limits_{0}^{2\pi} (u(t)-u(y))K(t-y) dy, \;\;\; t\in(0,2\pi),
 \end{equation}
and
 \begin{equation*}
   K(z)= \sum\limits_{n=- \infty}^{\infty}\dfrac{1}{| z - 2\pi n |^{1 + 2s}}, \;\;\; 0 < |z| <2\pi.
 \end{equation*}
 We consider the space $X$ defined in \cite{quaas} as the closure of the set of $2\pi$-periodic functions $u\in C^{1}(\mathbb{R})$
 with the norm
 \begin{equation*}
   \| u\|_{X} := \left( \dfrac{1}{2}\int\limits_{0}^{2\pi} \int\limits_{0}^{2\pi} (u(t)-u(y))^{2}K(t-y) dydt + \int\limits_{0}^{2\pi} u^{2}(t)dt\right)^{\dfrac{1}{2}},
 \end{equation*}
 where $X$ is a Hilbert space when provided with the inner product
 \begin{equation*}
   <u, v>:=  \dfrac{1}{2}\int\limits_{0}^{2\pi} \int\limits_{0}^{2\pi} (u(t)-u(y))(v(t) - v(y))K(t-y) dydt + \int\limits_{0}^{2\pi} u(t)v(t)dt.
 \end{equation*}
The space $X$ possesses good embedding properties which follow directly from the trivial relation $\|u\|_{H^{s}(0,2\pi)} \leq \|u\|_{X}$ for $u\in X$.\\

Consider now the family of equations
\begin{equation}\label{lanbda}
  (\triangle)^{s}v(t) + \lambda f(v(t))v'(t)=\lambda w(t) \quad \lambda \in (0,1).
\end{equation}
We want to establish Schauder type estimates for equation (\ref{lanbda}), but for this we need  first a $C^\alpha$  estimate that will follows from a $H^{s}$ bound and the embedding in  $C^\alpha$ by the fact that $s>1/2$.  For that, let us write
\begin{eqnarray*}
  v= \bar{v} + u \qquad&& \qquad w=\bar{w} + e,
\end{eqnarray*}
which implies that $\bar{u}=0$ and $\bar{e}=0$.  hence substituting this into (\ref{lanbda}) get
\begin{equation}\label{vecu}
(\triangle)^{s}u(t) +\lambda f(C + u(t))u'(t)=  \lambda e(t),
\end{equation}
  where $C= \bar{v}$.
  With that change we have the following Lemma.

\begin{lem}\label{lema21}
Let $u\in C^{2s + \alpha}$, for some $\alpha\in (0,1)$, be a classical solution  of (\ref{vecu}) with $\bar{e}=0$ and $\bar{u}=0$  then there exists $K>0$ such that
\begin{equation}
  \|u\|_{C^{\alpha}[0, 2\pi]} \leq K\|e\|_{L^{2}([0,2\pi])},
\end{equation}
for some $K>0$.
\end{lem}
\begin{dem}
Let $u$ is a $2\pi$-periodic solution of (\ref{vecu}). Multiplying by $u$ and integrating the equation (\ref{vecu}) we find
\begin{equation}\label{inter}
   \int\limits_{0}^{2\pi} (-\triangle)^{s}u u - \lambda \int\limits_{0}^{2\pi}  f(C + u)u' u= - \lambda \int\limits_{0}^{2\pi} eu.
\end{equation}
Using (\ref{lui}), together with the $2\pi$-periodicity of $u$ we find
\begin{equation*}
 [u]_{H^{s}(0,2\pi)}^{2}\leq [u]_{X}^{2} = \int\limits_{0}^{2\pi} u\mathcal{L} u =\int\limits_{0}^{2\pi} (-\triangle)^{s}u u
\end{equation*}

Let $H\in C^{2}(\mathbb{R})$ such that
\begin{equation*}
  H'(x)= \int\limits_{0}^{x} f(t)dt.
\end{equation*}
Then
\begin{eqnarray*}
 \int\limits_{0}^{2\pi} f(C + u(x))u'(t) u(t) dx &=& \int\limits_{0}^{2\pi} \dfrac{d}{dt}( H' (C + u(t)))u(t) dt \\
   &=& \int\limits_{0}^{2\pi} \frac{d}{dt}(H'(C + u(t))u(t))dt - \int\limits_{0}^{2\pi} H'(C + u(t))u'(t)dt\\
   &=& -\int\limits_{0}^{2\pi} \frac{d}{dt} H(C + u(t))dt=0.
\end{eqnarray*}

Hence substituting in (\ref{inter}) and  using Cauchy inequality we have
\begin{eqnarray}\label{abc}
  [u]_{H^{s}(0,2\pi)}^{2} &\leq&  \lambda\int\limits_{0}^{2\pi} |e(t)u(t)|dt \nonumber \\
               &\leq& \lambda \|e\|_{L^{2}}\|u\|_{L^{2}}.
\end{eqnarray}
On the other hand,  we claim that we have the following Poincaré type inequality. There exist $K>0$ such that
\begin{equation}\label{poincare}
   \int\limits_{0}^{2\pi}|u(t)|^{2}dt  \leq  K[u]_{H^{s}(0,2\pi)}^{2}.
\end{equation}
 Indeed, as $\bar{u}=0$ and  using Cauchy inequality  we get
\begin{eqnarray*}
  2\pi |u(t)| &=& \left| \int\limits_{0}^{2\pi} (u(t) -u(y))dy\right|  \\
              &\leq& \int\limits_{0}^{2\pi} |u(t) -u(y)|dy\\
  &\leq& \left(\int\limits_{0}^{2\pi} |u(t) -u(y)|^{2}dy\right)^{\dfrac{1}{2}}\left(\int\limits_{0}^{2\pi}dy\right)^{\dfrac{1}{2}}
\end{eqnarray*}

\begin{eqnarray*}
  |u(t)|^{2} &\leq& \dfrac{1}{2\pi}\int\limits_{0}^{2\pi} |u(t) -u(y)|^{2}dy  \\
               &=& \dfrac{1}{2\pi}\int\limits_{0}^{2\pi} \dfrac{|u(t) -u(y)|^{2}}{|t - y|^{1 + 2s}}|t - y|^{1 + 2s}dy,
\end{eqnarray*}
as $|t - y|< 2\pi$ gives
\begin{eqnarray*}
  \int\limits_{0}^{2\pi}|u(t)|^{2}dt &\leq& \dfrac{(2\pi)^{1 + 2s}}{2\pi}\int\limits_{0}^{2\pi}\int\limits_{0}^{2\pi}\dfrac{|u(t) -u(y)|^{2}}{|t - y|^{1 + 2s}}dydt \\
   &\leq& (2\pi)^{2s}[u]_{H^{s}(0,2\pi)}^{2},
\end{eqnarray*}
where $K= (2\pi)^{2s}$. Hence (\ref{abc}) and (\ref{poincare})  gives
\begin{equation*}
  \|u\|_{H^{s}(0,2\pi)} \leq \lambda K\|e\|_{L^{2}([0,2\pi])} \leq K\|e\|_{L^{2}([0,2\pi])},
\end{equation*}
and as $H^{s}(0,2\pi)\hookrightarrow C^{\alpha}[0,2\pi]$ for $s>\frac{1}{2}$ we get
\begin{equation}\label{desi}
  \|u\|_{C^{\alpha}[0,2\pi]} \leq K\|e\|_{L^{2}([0,2\pi])}.
\end{equation}
\fin
\end{dem}

 The following result is a basic property of our interest, where we will need the constant $C(1, s)$ which precisely given by
\begin{equation*}
  C(1,s)= \left( \int\limits_{\mathbb{R}} \dfrac{1 - \cos(\xi_{1})}{|\xi |^{1 + 2s}}d\xi\right)^{-1}.
\end{equation*}
 The proof can also be found in \cite{juan}.
\begin{lem}\label{perizero}
Let $u$ a $2\pi$-periodic function, then
\begin{equation*}
  \int\limits_{0}^{2\pi}  (-\triangle)^{s}u(x) dx= 0.
\end{equation*}
\end{lem}
\begin{dem}
Indeed, as $u$ is a periodic function, we can write its Fourier series:
\begin{equation*}
  u(x)= a_{0} + \sum\limits_{n=1}^{+ \infty} a_{n}\cos(nx) + b_{n}\sin(nx).
\end{equation*}
Then,
\begin{eqnarray}\label{cos}
  (-\triangle)^{s}u(x) &=& C_{1,s} \left(\sum\limits_{n=1}^{+ \infty} a_{n}\int\limits_{\mathbb{R}} \dfrac{\cos(nx) - \cos(ny)}{|x - y|^{1+ 2s}} dy + b_{n}\int\limits_{\mathbb{R}} \dfrac{\sin(nx) - \sin(ny)}{|x - y|^{1+ 2s}} dy \right)\nonumber\\
                       &=& C_{1,s} \left(\sum\limits_{n=1}^{+ \infty} a_{n}\int\limits_{\mathbb{R}} \dfrac{\cos(nx) - \cos(nx - nz)}{|z|^{1+ 2s}} dz + b_{n}\int\limits_{\mathbb{R}} \dfrac{\sin(nx) - \sin(nx -nz)}{|z|^{1+ 2s}} dz\right),\nonumber\\
                       &&
\end{eqnarray}
so
\begin{eqnarray*}
  \int\limits_{\mathbb{R}} \dfrac{\cos(nx) - \cos(nx - nz)}{|z|^{1+ 2s}} dz &=& \int\limits_{\mathbb{R}} \dfrac{\cos(nx) - \cos(nx)\cos(nz) - \sin(nx)\sin(nz)}{|z|^{1+ 2s}} dz \\
        &=& \int\limits_{\mathbb{R}} \dfrac{1 - \cos(nz)}{|z|^{1+ 2s}}dz\cos(nx) -  \int\limits_{\mathbb{R}}\dfrac{\sin(nz)}{|z|^{1 + 2s}}dz \sin(nx),
\end{eqnarray*}
where the second integral of the above is zero because we are integrating an odd function in a symmetric (with respect to zero) domain.
Similarly to the previous computation, we have
\begin{eqnarray*}
  \int\limits_{\mathbb{R}} \dfrac{\sin(nx) - \sin(nx - nz)}{|z|^{1+ 2s}} dz &=& \int\limits_{\mathbb{R}} \dfrac{\sin(nx) - \sin(nx)\cos(nz) + \cos(nx)\sin(nz)}{|z|^{1+ 2s}} dz \\
        &=& \int\limits_{\mathbb{R}} \dfrac{1 - \cos(nz)}{|z|^{1+ 2s}}dz\sin(nx) -  \int\limits_{\mathbb{R}} \dfrac{\sin(nz)}{|z|^{1 + 2s}}dz \cos(nx).
\end{eqnarray*}
Hence substituting this in  (\ref{cos}) gives
\begin{eqnarray*}
 (-\triangle)^{s}u(x) &=&  C_{1,s} \left(\sum\limits_{n=1}^{+ \infty} a_{n}\int\limits_{\mathbb{R}} \dfrac{1 - \cos(nz)}{|z|^{1+ 2s}}dz\cos(nx) + b_{n}\int\limits_{\mathbb{R}} \dfrac{1 - \cos(nz)}{|z|^{1+ 2s}}dz\sin(nx)\right) \\
   &=& \left( \int\limits_{\mathbb{R}} \dfrac{1 - \cos(\xi_{1})}{|\xi|^{n + 2s}}d \xi\right)^{-1} \left(\sum\limits_{n=1}^{+ \infty} a_{n}\int\limits_{\mathbb{R}} \dfrac{1 - \cos(nz)}{|z|^{1+ 2s}}dz\cos(nx)\right. \\
   && \qquad\qquad\qquad\qquad \qquad\qquad\qquad\qquad\left.+ b_{n}\int\limits_{\mathbb{R}} \dfrac{1 - \cos(nz)}{|z|^{1+ 2s}}dz\sin(nx)\right)\\
&=&  \sum\limits_{n=1}^{+ \infty} a_{n}n^{2s}\cos(nx) + b_{n}n^{2s}\sin(nx).
\end{eqnarray*}
Hence,
\begin{equation*}
   \int\limits_{0}^{2\pi} (-\triangle)^{s}u(x) dx = \sum\limits_{n=1}^{+ \infty} a_{n}n^{2s} \int\limits_{0}^{2\pi} \cos(nx)dx  + b_{n}n^{2s}  \int\limits_{0}^{2\pi} \sin(nx)dx= 0.
\end{equation*}
\fin
\end{dem}

\begin{dep}
Assume that $u$ is a $2\pi$-periodic solution of (\ref{vecu}), then integrating (\ref{vecu}) from $0$ to $2\pi$
\begin{equation*}
   \int\limits_{0}^{2\pi} (\triangle)^{s}u  +\lambda \int\limits_{0}^{2\pi}  f(C + u)u' =\lambda \int\limits_{0}^{2\pi} e,
\end{equation*}
hence, using the $2\pi$- periodicity  of $u$  and the notation of Lemma \ref{lema21} we have
\begin{equation*}
  \int\limits_{0}^{2\pi}  f(C+ u(t))u'(t)dt =  \int\limits_{0}^{2\pi} \dfrac{d}{dt} H' (C + u(t))dt =0,
\end{equation*}
which, together   with Lemma \ref{perizero} implies $\bar{e}=0$.

\medskip
For the existence, we will look at a priori estimates  for $u\in C^{1,\alpha}$, with $\alpha\in(0,1)$.
We first notice that  $f(C +u) \in C^{\alpha}(\mathbb{R})$, indeed as $u'$ is bounded we get $u$ is Lipschitz so
\begin{eqnarray}\label{fc}
 |f(C + u(t)) - f(C + u(y)|  &\leq& [f]_{\alpha}|u(t) - u(y)|^{\alpha} \nonumber\\
                               &\leq& [f]_{\alpha}|t - y|^{\alpha}\|u'\|_{C(\mathbb{R})}.
\end{eqnarray}
Now, $f(C +u)u' \in C^{\alpha}(\mathbb{R})$, because $f$ and $u'$ are bounded, using (\ref{fc}) we get
\begin{eqnarray}\label{fuc}
  |f(C + u(t))u'(t) - f(C + u(y))u'(y)| &\leq& |f(C + u(t)) - f(C + u(y)||u'(t)|\nonumber\\
  && + |f(C+ u(y))||u'(t) - u'(y)|\nonumber\\
    &\leq& [f]_{\alpha}|t -y|^{\alpha}\|u'\|_{C(\mathbb{R})}^{2} + \|f\|_{C(\mathbb{R})}|t -y|^{\alpha}[u']_{\alpha}.\nonumber\\
    &&
  \end{eqnarray}
Hence (\ref{fuc}) gives
\begin{eqnarray}\label{continua}
  \| f(C+u)u'\|_{C^{\alpha}(\mathbb{R})}   &\leq& \|f\|_{C^{\alpha}(\mathbb{R})} \|u\|_{C^{1,\alpha}(\mathbb{R})} \nonumber\\
   && + \|f\|_{C^{\alpha}(\mathbb{R})} \|u\|_{C^{1,\alpha}(\mathbb{R})}^{2} + \|f\|_{C^{\alpha}(\mathbb{R})}\|u\|_{C^{1,\alpha}(\mathbb{R})}\nonumber\\
&\leq& C\|u\|_{C^{1,\alpha}(\mathbb{R})}^{2} + C\|u\|_{C^{1,\alpha}(\mathbb{R})},
\end{eqnarray}
where $C=C(\|f\|_{C^{\alpha}(\mathbb{R})}, \alpha)$. As $s\in(\frac{1}{2},1)$ by [\cite{silvestre}, Proposition 2.9], we obtain that
\begin{eqnarray}\label{cotac1}
  \|u\|_{C^{1,\alpha}(\mathbb{R})} &\leq& C (\|u\|_{C(\mathbb{R})} + \|(\triangle)^{s}u\|_{C(\mathbb{R})})\nonumber \\
   &\leq& C(\|u\|_{C(\mathbb{R})}  +  \|\lambda f(C + u)u'\|_{C(\mathbb{R})} +  \|\lambda e\|_{C(\mathbb{R})})\nonumber\\
   &\leq&  C(\|u\|_{C(\mathbb{R})}  + \lambda \|f\|_{C(\mathbb{R})} \|u\|_{C^{1}(\mathbb{R})} +\lambda \|e\|_{C(\mathbb{R})})\nonumber\\
   &\leq& C(\|u\|_{C(\mathbb{R})}  + \|u\|_{C^{1}(\mathbb{R})} +\lambda \|e\|_{C(\mathbb{R})}),
\end{eqnarray}
where $C=C(\|f\|_{C(\mathbb{R})}, \alpha, s)$,  and using the Interpolation inequalities [\cite{interpol}, Theorems 3.2.1] we obtain that, for any $\epsilon > 0$, there is a positive constant $C = C(\alpha, \|f\|_{C^{\alpha}(\mathbb{R})},  \epsilon, s)$, such that
\begin{equation}\label{interc1}
  \|u\|_{C^{1}(\mathbb{R})}\leq \epsilon \|u\|_{C^{1, \alpha}(\mathbb{R})} + C\|u\|_{C(\mathbb{R})}.
\end{equation}
Hence, substituting (\ref{interc1}) in (\ref{cotac1}) and choosing $\epsilon= \frac{1}{2}$
\begin{equation}\label{c1c}
\|u\|_{C^{1,\alpha}(\mathbb{R})} \leq  C\|u\|_{C(\mathbb{R})} + \lambda C\|e\|_{C(\mathbb{R})}.
\end{equation}
Now, as $u$ is $2\pi$-periodic we have
\begin{equation*}
  \|u\|_{C(\mathbb{R})}= \|u\|_{C[0,2\pi]},
\end{equation*}
so that (\ref{desi}) implies
\begin{equation*}
  \|u \|_{C(\mathbb{R})}\leq \lambda K\|e\|_{L^{2}[0,2\pi]},
\end{equation*}
substituting this into (\ref{c1c}) we get
\begin{eqnarray}\label{co1}
  \|u \|_{C^{1,\alpha}(\mathbb{R})} &\leq& \lambda C(\|e\|_{C^{\alpha}(\mathbb{R})} + \|e\|_{L^{2}[0,2\pi]}) \nonumber.\\
                                        &\leq& \lambda C\|e\|_{C^{\alpha}(\mathbb{R})}.
\end{eqnarray}
Then, as $s>\frac{1}{2}$ by [\cite{silvestre}, Proposition 2.8],  we obtain that
\begin{equation}\label{estimate}
  \|u \|_{C^{2s + \alpha}(\mathbb{R})}\leq C(\| (\triangle)^{s} u\|_{C^{\alpha}(\mathbb{R})} +  \| u\|_{C(\mathbb{R})}),
\end{equation}
where $C=C(\alpha,s)$ is a positive constant. As $e\in C^{\alpha}(\mathbb{R})$ together with  (\ref{continua}) and (\ref{co1}) given
\begin{eqnarray}\label{estim}
   \|u \|_{C^{2s + \alpha}(\mathbb{R})} &\leq& C(  \|\lambda f(C+u)u'\|_{C^{\alpha}(\mathbb{R})} + \|\lambda e\|_{C^{\alpha}(\mathbb{R})} +  \|u\|_{C(\mathbb{R})})\nonumber \\
                                          &\leq& C(\lambda\|e\|_{C(\mathbb{R})}^{2} +\lambda\|e\|_{C^{\alpha}(\mathbb{R})} +  \| u\|_{C(\mathbb{R})}) \\
                                           &\leq& C(\|e\|_{C^{\alpha}(\mathbb{R})}^{2} +\|e\|_{C^{\alpha}(\mathbb{R})})= R,\nonumber
\end{eqnarray}
in this case $C= C(\alpha, \|f\|_{C^{\alpha}(\mathbb{R})}, s)$.\\

Now we complete the proof, if $z\in  C^{\alpha}_{2\pi} (\mathbb{R})$, the equation
\begin{equation*}
   (\triangle)^{s}u(t) + u(t)= z(t)
\end{equation*}
has a unique bounded classical solution $2\pi$ -periodic solution $u\in C^{2s + \alpha}(\mathbb{R})$,(see Lemma $3.1$, \cite{quaas}).\\
So we can define the map $K: C^{\alpha}_{2\pi} (\mathbb{R}) \to  C^{2s +\alpha}_{2\pi} (\mathbb{R})$ by $K(z)= u$ where $u$ is solution of
\begin{equation*}
   (\triangle)^{s}u(t) + u(t)= z(t),
\end{equation*}
so that $K: C^{\alpha}_{2\pi} (\mathbb{R}) \to  C^{1,\alpha}_{2\pi} (\mathbb{R})$ is compact because the injection of $C^{2s + \alpha}$ into $C^{1,\alpha}$ is compact since we have $s>1/2$. As a consequence, we defined $T: C^{1,\alpha}_{2\pi} (\mathbb{R}) \to  C^{1,\alpha}_{2\pi} (\mathbb{R})$ by $T(z)=u$  where $u$ is solution of
\begin{equation*}
  (\triangle)^{s}u(t) + u(t) = e(t) + z(t) - f(C+ z(t))z'(t)  ,
\end{equation*}
so $T$ is compact and using the Schaeffer theorem, $T$ has at least one fixed point $u$  in $B[0,R]\subset  C^{1,\alpha}_{2\pi} (\mathbb{R})$ so
\begin{equation*}
  u=T(u)
\end{equation*}
which implies that
\begin{equation*}
  (\triangle)^{s}u(t) =  e(t) - f(C+ u(t))u'(t).
\end{equation*}
 Consequently, $v= \bar{v} + u$ is a $2\pi$-periodic solution of (\ref{le}).
\fin
\end{dep}\\

Now we extend the previous results for systems of equations of the form
\begin{equation}\label{sistemc}
  L(u(t)) + \lambda \frac{d}{dt}( \nabla H(u(t))) + \lambda Au(t)= \lambda e(t).
\end{equation}
For simplicity to find a priori bounds of the system $(\ref{sistemc})$ only in the case $\lambda=1$.\\

\begin{desis}
 First, we reduce the case a solutions with zero mean value, let us write
\begin{equation*}
  \begin{array}{cc}
  u(t)= \bar{u} + v(t), & e(t)= \bar{e} + w(t)
\end{array}
\end{equation*}
which implies that $\bar{v}=0$ and $\bar{e}=0$, then substituting this into (\ref{sistema}) becomes
\begin{equation}\label{ecuaresol}
  L(v(t)) + \frac{d}{dt}( \nabla  H(\bar{u} +v(t)))+ A\bar{u} + Av(t)= \bar{e} + w(t),
\end{equation}
we get the equivalent system
\begin{equation}\label{bar}
  A\bar{u}= \bar{e},
\end{equation}
\begin{equation}\label{ecuafin}
   L(v(t)) + \frac{d}{dt}( \nabla H(\bar{u} +v(t))) + Av(t)=  w(t).
\end{equation}
As $\bar{e}\in Im A$, then exist at least one $\bar{u}$ solving (\ref{bar}), so for each $\bar{u}$, we just have to find a $2\pi$-periodic solution of (\ref{ecuafin}) $v\in C^{1, \alpha}( \mathbb{R}, \mathbb{R}^{n})$.\\

On the other hand, note that
\begin{equation*}
  \frac{d}{dt}( \nabla H(\bar{u} +v(t))) = F(\bar{u} + v(t))v'(t),
\end{equation*}
 where $F$ is the Hessian matrix of H, so substituting in (\ref{ecuafin}) get
\begin{equation}\label{sol}
   L(v(t)) + F(\bar{u} + v(t))v'(t) + Av(t)=  w(t).
\end{equation}
Now, we shall find a priori estimates, first will look estimates for $H^{s}((0, 2\pi),\mathbb{R}^{n})$ where $s=\min\limits_{1\leq i\leq n} s_{i}$. We will use the usual norm, but it is important to mention that for $s \leq s_{i}$ we have
\begin{eqnarray}\label{hss}
  \|v\|_{H^{s}}^{2} &=&  \sum\limits_{1}^{n}   \|v_{i}\|_{H^{s}}^{2} \\
   &\leq& \sum\limits_{1}^{n}   \|v_{i}\|_{H^{s_{i}}}^{2}.
\end{eqnarray}
 Now, assume $v$ is a $2\pi$-periodic solution of (\ref{sol}), multiplying for $v$ and integrating the equation
\begin{equation*}
 \int\limits_{0}^{2\pi} <-L(v(t)),v(t)>dt - \int\limits_{0}^{2\pi} <F(\bar{u} + v(t))v'(t),v(t)> dt - \int\limits_{0}^{2\pi} <Av(t),v(t)>dt
\end{equation*}
\begin{equation}\label{sistemul}
  = -\int\limits_{0}^{2\pi} <w(t),v(t)>dt,
\end{equation}
so that $(\ref{hss})$ together the $2\pi$-periodicity of $v$ gives
\begin{equation*}
 \int\limits_{0}^{2\pi}  <-L(v(t)), v(t)> =  \int\limits_{0}^{2\pi} \sum\limits_{1}^{n}  (-\triangle) ^{s_{i}}v_{i}(t)v_{i}(t),
\end{equation*}
so
\begin{eqnarray}\label{enrn}
 [v]_{H^{s}((0, 2\pi),\mathbb{R}^{n})}^{2} &\leq&  \sum\limits_{1}^{n}  [v_{i}]_{H^{s_{i}}((0, 2\pi)}^{2}\\
  &\leq&  \sum\limits_{1}^{n}  \int\limits_{0}^{2\pi} (-\triangle) ^{s_{i}}v_{i}(t)v_{i}(t)\\
  &=&  \int\limits_{0}^{2\pi}  <-L(v(t)), v(t)>.
\end{eqnarray}
\begin{eqnarray*}
 \int\limits_{0}^{2\pi} <F(\bar{u} + v(t))v'(t),v(t)> &=& \int\limits_{0}^{2\pi} <\frac{d}{dt}( \nabla H(v(t))),v(t)>dt\\
                                               &=& \int\limits_{0}^{2\pi} \dfrac{d}{dt} <  \nabla  H(v(t)),v(t)> dt\\
                                                 && - \int\limits_{0}^{2\pi}  <  \nabla  H(v(t)),v'(t)> dt\\
                                                &=& \int\limits_{0}^{2\pi} \frac{d}{dt} H(v(t))dt=0.
\end{eqnarray*}
 Hence substituting in (\ref{sistemul}) gives
\begin{equation}\label{hol}
  [v]_{H^{s}((0,2\pi),\mathbb{R}^{n})}^{2} -  \int\limits_{0}^{2\pi} <Av(t),v(t)>dt \leq - \int\limits_{0}^{2\pi} < w(t),v(t)>dt.
\end{equation}
By assumption we have
\begin{equation*}
   \int\limits_{0}^{2\pi} <Av(t),v(t)>dt \leq M \| v\|_{L^{2}}^{2},
\end{equation*}
so that (\ref{hol}) implies using H\"older inequality implies
\begin{equation}\label{desial}
    [v]_{H^{s}((0,2\pi),\mathbb{R}^{n})}^{2} \leq M \| v\|_{L^{2}}^{2} + \| w\|_{L^{2}}\| v\|_{L^{2}}.
\end{equation}
Therefore, the Poincaré type inequality gives
\begin{equation*}
   [v]_{H^{s}((0,2\pi),\mathbb{R}^{n})} \leq (1 -M)^{-1}\| w\|_{L^{2}},
\end{equation*}
and hence
\begin{equation}\label{estimax}
  \|v\|_{H^{s}((0,2\pi),\mathbb{R}^{n})} \leq K \| w\|_{L^{2}}.
\end{equation}

Then, we shall find a priori estimates in H\"older spaces; as $s>\frac{1}{2}$ by [\cite{silvestre}, Proposition 2.8] gives
\begin{equation}\label{cotaensist}
  \|v_{i} \|_{C^{2s + \alpha}(\mathbb{R})}\leq C(\| (\triangle)^{s} v_{i}\|_{C^{\alpha}(\mathbb{R})} +  \| v_{i}\|_{C(\mathbb{R})}),
\end{equation}\label{2cotasist}
hence we get
\begin{eqnarray}\label{cotaprinc}
  \| (\triangle)^{s_{i}} v_{i}\|_{C^{\alpha}(\mathbb{R})} &\leq& C\sum\limits_{i=1}^{n}\|v_{i}\|_{C^{1,\alpha}(\mathbb{R})}^{2} + C\sum\limits_{i=1}^{n}\|v_{i}\|_{C^{\alpha}(\mathbb{R})} +  C\|w\|_{C^{\alpha}(\mathbb{R})} \nonumber \\
  &\leq& C(\sum\limits_{i=1}^{n}\|v_{i}\|_{C^{1,\alpha}(\mathbb{R})}^{2} + \sum\limits_{i=1}^{n}\|v_{i}\|_{C^{1,\alpha}(\mathbb{R})} +  \|w\|_{C^{\alpha}(\mathbb{R})}),
\end{eqnarray}
where $C=C(n,s, \alpha, \|f\|_{C^{\alpha}(\mathbb{R})}, \|A\|)$, by [\cite{silvestre}, Proposition 2.9] we obtain
\begin{eqnarray}\label{29}
    \| v_{i}\|_{C^{1,\alpha}(\mathbb{R})} &\leq& C(\sum\limits_{i=1}^{n}\|v_{i}\|_{C(\mathbb{R})} + \| (\triangle)^{s} v_{i}\|_{C(\mathbb{R})})\\
      &\leq& C(\sum\limits_{i=1}^{n}\|v_{i}\|_{C(\mathbb{R})} + \sum\limits_{i=1}^{n}\|v_{i}\|_{C^{1}(\mathbb{R})}+ \|w\|_{C(\mathbb{R})}),
\end{eqnarray}
 and using the Interpolation inequalities [\cite{interpol}, Theorems 3.2.1] together with the fact that $2s > 1$, we obtain that, for any $\epsilon > 0$, there is a positive constant $C = C(n,\alpha, \|f\|_{C^{\beta}(\mathbb{R})}, \|A\|, \epsilon, s)$, such that
\begin{eqnarray*}
    \| v_{i}\|_{C^{1}(\mathbb{R})} &\leq& \epsilon \sum\limits_{i=1}^{n} \|v_{i}\|_{C^{1,\alpha}(\mathbb{R})} + C\sum\limits_{i=1}^{n}\|v_{i}\|_{C(\mathbb{R})}.
\end{eqnarray*}
Hence,
\begin{equation}\label{sum}
 \sum\limits_{i=1}^{n} \| v_{i}\|_{C^{1}(\mathbb{R})} \leq \epsilon \sum\limits_{i=1}^{n}\sum\limits_{i=1}^{n} \|v_{i}\|_{C^{1,\alpha}(\mathbb{R})} + C\sum\limits_{i=1}^{n}\sum\limits_{i=1}^{n}\|v_{i}\|_{C(\mathbb{R})}.
\end{equation}
Thus,
\begin{eqnarray*}
   \sum\limits_{i=1}^{n} \| v_{i}\|_{C^{1,\alpha}(\mathbb{R})} &\leq& C\sum\limits_{i=1}^{n}\sum\limits_{i=1}^{n}\|v_{i}\|_{C(\mathbb{R})} + 2n\epsilon C \sum\limits_{i=1}^{n} \|v_{i}\|_{C^{1,\alpha}(\mathbb{R})} + C\|w\|_{C(\mathbb{R})}.
\end{eqnarray*}

Now, as $v$ is $2\pi$-periodic we have
\begin{equation*}
  \|v_{i}\|_{C(\mathbb{R})}= \|v_{i}\|_{C[0,2\pi]},
\end{equation*}
so that (\ref{desi}) implies
\begin{equation*}
  \|v_{i} \|_{C(\mathbb{R})} \leq  K\|w\|_{L^{2}(\mathbb{R})}.
\end{equation*}
 Hence choosing $\epsilon = \dfrac{1}{4nC}$, we get
\begin{eqnarray}\label{3cota}
\sum\limits_{i=1}^{n} \| v_{i}\|_{C^{1,\alpha}(\mathbb{R})} &\leq& C\|w\|_{L^{2}(\mathbb{R})} + C\|w\|_{C(\mathbb{R})}\\
                                                           &\leq& C \|w\|_{C(\mathbb{R})}.
\end{eqnarray}
Now,
\begin{equation}\label{cotacua}
  \sum\limits_{i=1}^{n} \| v_{i}\|_{C^{1,\alpha}(\mathbb{R})}^{2} \leq \left(\sum\limits_{i=1}^{n} \| v_{i}\|_{C^{1,\alpha}(\mathbb{R})}\right)^{2} \leq \|w\|_{C(\mathbb{R})}^{2}.
\end{equation}
Hence substituting (\ref{3cota}) and (\ref{cotacua}) into (\ref{cotaprinc}) gives
\begin{equation}\label{pq}
   \| (\triangle)^{s_{i}} v_{i}\|_{C^{\alpha}(\mathbb{R})} \leq C\|w\|_{C(\mathbb{R})}^{2} + C\|w\|_{C(\mathbb{R})}.
\end{equation}
Substituting (\ref{pq}) in (\ref{cotaensist}) we get
\begin{eqnarray*}
   \| v_{i}\|_{C^{2s + \alpha}(\mathbb{R})}&\leq& C\|w\|_{C(\mathbb{R})}^{2} + C\|w\|_{C(\mathbb{R})}+ \|v_{i}\|_{C(\mathbb{R})}\\
   &\leq& C\|w\|_{C(\mathbb{R})}^{2} + C\|w\|_{C(\mathbb{R})},
\end{eqnarray*}
so
\begin{eqnarray*}
  \sum\limits_{i=1}^{n} \| v_{i}\|_{C^{2s + \alpha}(\mathbb{R})}&\leq& C\|w\|_{C^{\alpha}(\mathbb{R})}^{2} + C\|w\|_{C(\mathbb{R})}= R,
\end{eqnarray*}
where $C=C(n, s, \alpha, \|f\|_{C^{\alpha}(\mathbb{R})}, \|A\|)$. Now, if $y_{i}\in C^{1, \alpha}_{2\pi}(\mathbb{R})$ the equation
\begin{equation*}
   (\triangle)^{s_{i}}v_{i}(t) + v_{i}(t)= y_{i}(t) + v_{i}(t)
\end{equation*}
has a unique bounded classical solution $2\pi$ -periodic solution $v_{i}\in C^{2s + \alpha}_{2\pi}(\mathbb{R})$, then
exist $v(t)=(v_{1}(t), \cdots,  v_{n}(t))$ such that is solution of
\begin{equation*}
  L(v(t)) + v(t)= y(t) + v(t),
\end{equation*}
where $y(t)= (y_{1}(t),\cdots,  y_{n}(t))$ .\\
Finally, we define the operator  $T: C^{1, \alpha}_{2\pi}(\mathbb{R}, \mathbb{R}^{n}) \to C^{1, \alpha}_{2\pi}(\mathbb{R}, \mathbb{R}^{n})$ by $T(y)=v$ where $v$ is solution the
\begin{equation*}
  L(v(t)) + v(t)=w(t) + y(t) - F(\bar{u} + y(t))y'(t) - Ay(t),
\end{equation*}
 this operator is compact, because the injection of $C^{2s + \alpha}_{2\pi}(\mathbb{R}, \mathbb{R}^{n})$ into  $C^{1, \alpha}_{2\pi}(\mathbb{R}, \mathbb{R}^{n})$ is compact, so using the Schaeffer theorem, $T$ has at least one fixed point $v$ in $B(0,R) \subset C^{1, \alpha}_{2\pi}(\mathbb{R}, \mathbb{R}^{n})$, this  is
 \begin{equation*}
   T(v)=v
 \end{equation*}
which implies that
\begin{equation*}
  L(v(t))= w(t) - F(\bar{u} + v(t))v'(t) - Av(t),
\end{equation*}
i.e. $v$ is a $2\pi$-periodic solution of (\ref{sol}). Consequently, $u= \bar{u} + v$ is a $2\pi$-periodic solution of (\ref{sistema}).
\fin
\end{desis}

\section{Case of attractive singularity}

The standard method of sub and super solutions provides the following existence theorem for the $2\pi$-periodic solutions of equation (\ref{intro}), the more difficult part is to construct super-solution and is the place where we use Proposition \ref{l1}.

\begin{lem}\label{l2}
Assume that there exist $2\pi$- periodic $C^{2s + \alpha}$- functions  $\eta$ and $\beta$ for some $\alpha\in(0,1)$ and $s\in(1/2,1)$ such that $\eta\leq\beta$ and
\begin{eqnarray*}
(\triangle)^{s}\eta + f(\eta)\eta' +g(\eta) &\geq& e \\
   (\triangle)^{s}\beta + f(\beta)\beta' +g(\beta) &\leq& e,
\end{eqnarray*}
in $\mathbb{R}$. Then equation (\ref{intro}) has at least one $2\pi$-periodic classical solution $u\in C^{2s +\alpha}_{2\pi}(\mathbb{R})$ satisfying  $\eta(x)\leq u(x) \leq\beta(x)$ for all $x\in \mathbb{R}$.
\end{lem}

\begin{dem}

Let $u\in C^{1,\alpha}_{2\pi}(\mathbb{R})$, we define for each $x\in \mathbb{R}$
\begin{equation*}
  H(u(x),u'(x))= \begin{cases}
-\beta(x) + e(x) - f(\beta(x))\beta'(x) - g(\beta(x)) , & \text{if $u(x)>\beta(x)$},\\
-u(x) + e(x) - f(u(x))u'(x) -  g(u(x)), & \text{if $\eta(x)\leq u(x)\leq \beta(x)$},\\
-\eta (x) + e(x) - f(\eta(x))\eta'(x) -g(\eta(x)), & \text{if $u(x)<\eta(x)$,}
\end{cases}
\end{equation*}
then $H\in C^{\alpha}_{2\pi}(\mathbb{R})$ .\\

Thus we can define
\begin{eqnarray*}
  N: C^{1,\alpha}_{2\pi} (\mathbb{R})&\to& C^{\alpha}_{2\pi}(\mathbb{R})  \\
  u &\longmapsto& H(u,u')
\end{eqnarray*}
and
\begin{eqnarray*}
  K: C^{\alpha}_{2\pi} (\mathbb{R})&\to& C^{1,\alpha}_{2\pi}(\mathbb{R})  \\
              z &\longmapsto& v
\end{eqnarray*}
where $v$ is the unique solution of
\begin{equation*}
  (\triangle)^{s}v(x) - v(x)=z(x).
\end{equation*}
Then, we have $K \circ N: C^{1,\alpha}_{2\pi} (\mathbb{R})\to C^{1,\alpha}_{2\pi} (\mathbb{R})$ is continuous and compact, as proved in Proposition (\ref{l1}).
Notice that $N(C^{1,\alpha}_{2\pi} (\mathbb{R}))$ is bounded. Hence, by Schauder's fixed point theorem, $K\circ N$ has a fixed point $u$, i.e., $u$ is a solution of
\begin{equation*}
  (\triangle)^{s}u(x) - u(x)= H(u,u') \;\;\; in \;\;\; \mathbb{R}.
\end{equation*}
Now we prove that $\eta(x)\leq u(x)\leq \beta(x)$ for all $x\in \mathbb{R}$, we first show that $u(x)\leq \beta(x)$. The other inequality is similar.
We assume by contradiction that $\max (u(x) - \beta(x))= u(\bar t) - \beta(\bar t)> 0$, here $\bar t$ is the point where the maximum is attained.
So we have
 $(\triangle)^{s}u(\bar t) - (\triangle)^{s}\beta (\bar t) \leq 0$ and also
\begin{eqnarray*}
 (\triangle)^{s}u(\bar t) -(\triangle)^{s}\beta (\bar t)  &\geq & u(\bar t) + H(u(\bar t))+ f(\beta(\bar t))\beta'(\bar t) + g(\beta(\bar t)) - e(\bar t)\\
 &\geq& u(\bar t) -\beta(\bar t) + e(\bar t)- f(\beta(\bar t))\beta'(\bar t)- g(\beta(\bar t))+ f(\beta(\bar t))\beta'(\bar t)+ g(\beta(\bar t)) - e(\bar t)\\
 &=& u(\bar t) - \beta(\bar t)>0,
\end{eqnarray*}
this gives a contradicts . Therefore, $u(x)\leq \beta(x)$.
\fin
\end{dem}\\

 Now we are in position to prove Theorem \ref{teo}, which is direct now.\\

\begin{deteo}
By assumption $1$, there exist a constant $\eta >0$ such that
\begin{equation*}
  g(\eta)\geq e(x) \quad \text{in} \quad [0,2\pi],
\end{equation*}
and thus $\eta$ is a sub solution for (\ref{intro}) with $2\pi$-periodic. We now write $e(x)= \bar{e} + \tilde{e}$, then by assumption $2$, there exists $R>0$ be such that $g(x) \leq \bar{e}$ for $x \geq R$ and by Proposition \ref{l1}, the equation
\begin{equation*}
 (\triangle)^{s}v(x) + f(C + v(x))v'(x)=\tilde{e}(x)
\end{equation*}
has one $2\pi$-periodic solution $v$, so we take $C$ sufficiently largo such that $C + v(x)\geq \max (\eta, R)$ for all $x\in [0, 2\pi]$. Hence we take $\beta(x)= C + v(x)$, gives
\begin{eqnarray*}
   (\triangle)^{s}\beta(x) + f(\beta)\beta' +g(\beta) &=&  (\triangle)^{s}v(x) + f(C + v(x))v'(x) +g(C + v(x)) \\
                                                      &\leq& \tilde{e}(x) + \bar{e}= e(x),
\end{eqnarray*}
so that $\beta(x)\geq \eta$ is a super-solution for (\ref{intro}) with $2\pi$-periodic. Then using Lemma \ref{l2} there exists a $2\pi$-periodic solution $u$ of (\ref{intro}) with $\eta \leq u(x) \leq \beta(x)$.
\fin
\end{deteo}
\begin{coro}\label{corot}
Assume that the function $g:(0, + \infty) \to (0, + \infty)$ is such that the following conditions hold.
\begin{enumerate}
  \item $g(x) \to + \infty$ as $x \to 0^{+}$,
  \item $\limsup\limits_{x  \to + \infty} g(x) = 0$.
\end{enumerate}
Then, if $\bar{e}>0$ the equation (\ref{intro}) has a positive $2\pi$- periodic solution.
\end{coro}
\begin{dem}
Follows directly from Theorem \ref{teo}.
\fin
\end{dem}\\

\begin{example} The Forbat type equation fractional Laplacian that is
\begin{equation}\label{qwe}
  (\triangle)^{s}v(x) + f(v(x))v'(x) + \dfrac{v(x)}{v(x) - C}=  e(x),
\end{equation}
where $f\in C^{\alpha}(\mathbb{R})$ and $e\in C^{\alpha}(\mathbb{R})$ and $\overline{e -1}>0$. By the change of variable $u=v - C$ we have
\begin{equation}\label{ex}
    (\triangle)^{s}u(x) + f(C + u(x))u'(x) + \dfrac{C}{u(x)}= e(x) -1,
\end{equation}
so let us now take
\begin{equation*}
  g(x)=\dfrac{C}{x},
\end{equation*}
the Corollary \ref{corot} implies that equation (\ref{ex}) has at least one $2\pi$-periodic positive solution, i.e. the equation (\ref{qwe}) has at least one $2\pi$-periodic solution $v$ such that $v (x)>C$ for all $x\in[0,2\pi]$.
\end{example}

\section{Case of repulsive singularity}

In this section, we will prove the existence of a positive $2\pi$-periodic solution to (\ref{ecuag2}). For that, let as consider the family of equations:
\begin{equation}\label{ecu2l}
   (-\triangle) ^{s}u(t) + cu'(t) - \lambda g(u) =\lambda e(t), \qquad \lambda \in (0,1).
\end{equation}

First, we want to find a priori bounds. For that the following Lemma is very important and corresponds  to a new energy type identity for the periodic solutions, which
proof is based in ideas of Lemma $3.2$ in \cite{simetria}.


\begin{lem}\label{lemint}
	Let $u\in C^{1}_{2\pi}(\R)$, then
	\begin{eqnarray}\label{ener}
	  \int\limits_{a}^{b} u'(x) (-\Delta)^s u(x) dx &=&\frac{c(1,s)}{2} \left( \int\limits_{-\infty}^{+\infty} \frac{(u(b)-u(y))^2}{|b-y|^{1+2s}} dy   - \int\limits_{-\infty}^{+\infty} \frac{(u(a)-u(y))^2}{|a -y|^{1+2s}} dy \right.\nonumber\\
	   &&  \left.  - (1+2s)\int\limits_{a}^{b} \int\limits_{b}^{+ \infty} \frac{(u(x)-u(y))^2}{|x-y|^{2+2s}}dydx\right.\nonumber\\
     &&  \left. + (1+2s)\int\limits_{a}^{b} \int\limits_{-\infty}^{a} \frac{(u(x)-u(y))^2}{|x-y|^{2+2s}} dy dx\right),
	\end{eqnarray}
	for every $0<a<b$.
	\end{lem}

	\begin{dem}
	Fix $0<a<b$ and choose $\delta$ and $M$ with the restrictions $0<\delta<a$ and $M>b+\delta$. We
	consider the integral
	\begin{equation}\label{Idelta}
	I_{\delta,M}= \int_{a}^{b} u'(x) \int^M \limits_{-\text{\em\scriptsize M} \atop {|y-x|\ge \delta}} \frac{u(x)-u(y)}{|x-y|^{1+2s}} dy dx
	= \iint_{A_{\delta,M}}u'(x) \frac{u(x)-u(y)}{|x-y|^{1+2s}} dy dx,
	\end{equation}
	where $A_{\delta,M}=([a,b]\times [-M,M]) \cap \{(x,y)\in \R^2:\ |y-x|\ge \delta\}$.

It is not hard to see that
	\begin{align*}
	I_{\delta,M} & = \frac{1}{2} \iint_{A_{\delta,M}} \hspace{-2mm} \left(\dfrac{(u(x)-u(y))^2}{|x-y|^{1+2s}}\right)_x
	\hspace{-1mm} dy dx + \frac{1+2s}{2} \hspace{-1mm}
	\iint_{A_{\delta,M}} \hspace{-3mm} \frac{(x-y) (u(x)-u(y))^2}{|x-y|^{3+2s}} dy dx.
	\end{align*}
	
    We now split $A_{\delta,M}=A^1 \cup A^2 \cup A^3 \cup A^4$, where
	$$
	\begin{array}{l}
	A^1= \{(x,y) \in A_{\delta,M}: y \ge b\}\\[0.25pc]
	A^2= \{(x,y) \in A_{\delta,M}: x + \delta \le y \le b\}\\[0.25pc]
    A^3= \{(x,y) \in A_{\delta,M}: a \le y \le x - \delta\}\\[0.25pc]
	A^4= \{(x,y) \in A_{\delta,M}: y\le a\}.
	\end{array}
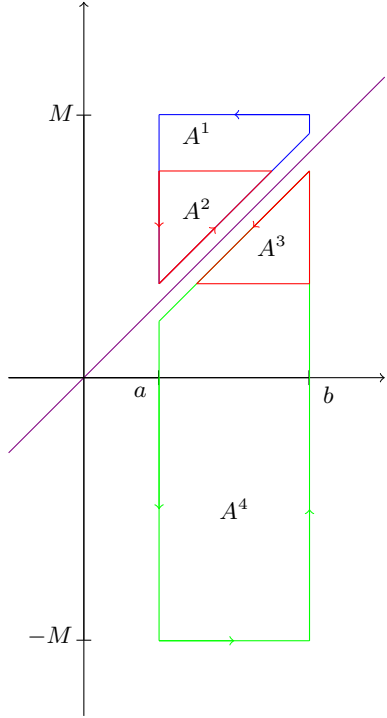
	$$
\begin{center}
	\begin{tikzpicture}

\node [below] at (0.75,0) {$a $};

\node [below] at (3.25,0) {$b$};

\node [below] at (-0.3,3.75) {$M $};

\node [below] at (-0.45,-3.2) {$-M$};

\node [below] at (2,-1.5) {$A^4$};

\node [below] at (2.5,2) {$A^3$};

\node [below] at (1.5,2.5) {$A^2$};

\node [below] at (1.5,3.5) {$A^1$};

\draw  [->]  (-1,0) -- (4,0)  ;

\draw   [->] (0,-4.5) -- (0,5)  ;

\draw  [violet]  (-1,-1) -- (4,4)  ;

\draw  [-|]  (-1,0) -- (1,0)  ;

\draw  [-|]  (-1,0) -- (3,0)  ;

\draw  [-|]  (0,0) -- (0,3.5)  ;

\draw  [-|]  (0,0) -- (0,-3.5) ;

\draw  [green]  (1,0.75) -- (3,2.75)  ;

\draw  [green]  (1,0.75) -- (1,-3.5)  ;

\draw  [green] [->] (1,0) -- (1,-1.75)  ;

\draw  [green]  (3,2.75) -- (3,-3.5)  ;

\draw  [green] [->]  (3,-3.5) --(3,-1.75)   ;

\draw  [green] [->]  (1,-3.5) --(2,-3.5)   ;

\draw  [green] (1,-3.5) --(3,-3.5)   ;

\draw  [blue]  (1,3.5) -- (3,3.5)  ;

\draw  [blue] [<-] (2,3.5)  --   (3,3.5) ;

\draw  [blue]  (1,3.5) -- (1,1.25)  ;

\draw  [blue]  (1,1.25) -- (3,3.25)  ;

\draw  [blue]  (3,3.5) -- (3,3.25)   ;

\draw  [red]   (1.5,1.25)  --   (3,1.25);

\draw  [red]    (3,1.25)   --  (3,2.75) ;

\draw  [red]   (1.5,1.25)  --  (3,2.75) ;

\draw  [red]  [<-]  (2.25,2)  --  (3,2.75) ;

\draw  [red]   (1,2.75)  --   (2.5,2.75);

\draw  [red]   (1,2.75)  --   (1,1.25);

\draw  [red]    (1,1.25) --(2.5,2.75);

\draw  [red]  [->]   (1,1.25) --(1.75,2);

\draw  [red] [->]   (1,2.75)  --   (1,2);

\end{tikzpicture}
	
	\end{center}

	Since the region $A^2$ is the reflection of $A^3$, with respect to the line $y=x$ and the integrand in the last integral above is antisymmetric, we get  	
	\begin{align*}
	I_{\delta,M} &= \frac{1}{2}\iint\limits_{A_{\delta,M}} \left(\frac{(u(x)-u(y))^2}{|x-y|^{1+2s}}\right)_x dydx - \frac{1+2s}{2}\left( \iint\limits_{A^1} \frac{(u(x)-u(y))^2}{(x-y)^{2+2s}} dy dx - \iint\limits_{A^4}\frac{(u(x)-u(y))^2}{(x-y)^{2+2s}} dy dx \right)\\
	& =\frac{1}{2}\oint\limits_{\partial A_{\delta,M}} \frac{(u(x)-u(y))^2}{|x-y|^{1+2s}}dy- \frac{1+2s}{2}\left( \iint\limits_{A^1} \frac{(u(x)-u(y))^2}{(x-y)^{2+2s}} dy dx - \iint\limits_{A^4}\frac{(u(x)-u(y))^2}{(x-y)^{2+2s}} dy dx \right).
	\end{align*}
	Now, using Green's formula, where the line integral is to be taken in the positive sense. Parameterizing
	the line integral we have

    \begin{eqnarray}\label{eq-energis-1}
     I_{\delta,M}  &=& -\frac{1}{2}\int^M \limits_{-\text{\em\scriptsize M}  \atop |y-a|\ge \delta} \frac{(u(a)-u(y))^2}{|a-y|^{1+2s}} dy
	+\frac{1}{2}\int^M \limits_{-\text{\em\scriptsize M}  \atop |y-b|\ge \delta} \frac{(u(b)-u(y))^2}{|b-y|^{1+2s}} dy \nonumber  \\
      && + \frac{1}{2}\int^{b-\delta}\limits_{a} \frac{(u(x)-u(x+\delta))^2}{\delta^{1+2s}} dx
	-\frac{1}{2}\int_{a}^{b} \frac{(u(x)-u(x-\delta))^2}{\delta^{1+2s}} dx \nonumber\\
       && - \frac{1+2s}{2} \iint_{A^1} \frac{(u(x)-u(y))^2}{(x-y)^{2+2s}} dy dx + \frac{1+2s}{2} \iint_{A^4} \frac{(u(x)-u(y))^2}{(x-y)^{2+2s}} dy dx\nonumber\\
       &=&  -\frac{1}{2}\int^M \limits_{-\text{\em\scriptsize M} \atop |y-a|\ge \delta} \frac{(u(a)-u(y))^2}{|a-y|^{1+2s}}dy +\frac{1}{2}\int^M \limits_{-\text{\em\scriptsize M}  \atop |y-b|\ge \delta} \frac{(u(b)-u(y))^2}{|b-y|^{1+2s}} dy  \nonumber \\
       &&  - \frac{1}{2}\int\limits_{a-\delta}^{a}  \frac{(u(x+\delta)-u(x))^2}{\delta^{1+2s}} dx
		- \frac{1+2s}{2} \iint_{A^1} \frac{(u(x)-u(y))^2}{(x-y)^{2+2s}} dy dx\nonumber\\
        && + \frac{1+2s}{2} \iint_{A^4} \frac{(u(x)-u(y))^2}{(x-y)^{2+2s}} dy dx.
    \end{eqnarray}

	
	Now we can pass to the limit in $I_{\delta, M}$ as $M\to +\infty$ using dominated convergence and get
	
	\begin{eqnarray}\label{eq-energia-2}
	\int_{a}^{b} u'(x) \int_{|y-x|\ge \delta} \frac{u(x)-u(y)}{|x-y|^{1+2s}} dy dx
	&=& \frac{1}{2}\int_{|y-b|\ge \delta} \frac{(u(b)-u(y))^2}{|b-y|^{1+2s}} dy\nonumber\\
 && -\frac{1}{2}\int_{|y-a|\ge \delta} \frac{(u(a)-u(y))^2}{|a-y|^{1+2s}} dy \nonumber \\
	&& - \frac{1}{2}\int_{a-\delta}^{a}  \frac{(u(x+\delta)-u(x))^2}{\delta^{1+2s}} dx \\
    && - \frac{1+2s}{2} \iint_{A_\delta^{1}} \frac{(u(x)-u(y))^2}{(x-y)^{2+2s}} dy dx,\nonumber\\
	&&+ \frac{1+2s}{2} \iint_{A_\delta^{2}} \frac{(u(x)-u(y))^2}{(x-y)^{2+2s}} dy dx,\nonumber
	\end{eqnarray}
	where $A_\delta^{1}=([a, b]\times (b, +\infty) \cap \{(x,y)\in \R^2: y \geq x + \delta\}$ and  $A_\delta^{2}=([a, b]\times (- \infty, a) \cap \{(x,y)\in \R^2: y \le x-\delta\}$.
	
	Finally, will want to pass to the limit as $\delta \to 0$ in \eqref{eq-energia-2}. Observe
	that, since $u\in C^1_{2\pi} (\R)$, we have that $u'$ is bounded so $u$ is Lipschitz, hence for $y$ close to $b$ gives
	$$
	\frac{(u(b)-u(y))^2}{|b-y|^{1+2s}} \le C |b-y|^{1-2s}\in L^1_{\rm loc}(\R),
	$$
and we also have for $y$ close to $a$
	$$
	\frac{(u(a)-u(y))^2}{|a-y|^{1+2s}} \le C |a-y|^{1-2s}\in L^1_{\rm loc}(\R).
	$$

	So the passing to the limit is justified in the first and second integral in the right-hand side of
	\eqref{eq-energia-2} by dominated convergence. As for the third integral, for being $u$ Lipschitz gives
	$$
    \frac{(u(x+\delta)-u(x))^2}{\delta^{1+2s}} \le C\delta^{1-2s}
    $$
    so that,
	$$
	\int_{t-\delta}^{t} \frac{(u(x+\delta)-u(x))^2}{\delta^{1+2s}} dx \le C\delta^{2-2s}\to 0
	$$
	as $\delta \to 0^+$. As for the double integral, we also have that
	$$
	\frac{(u(x)-u(y))^2}{|x-y|^{2+2s}} \le C |x-y|^{-2s}\in L^1_{\rm loc}(\R^2),
	$$
	for $x$ and $y$ close.  Therefore, we can pass to the limit in the
	right-hand side of \eqref{eq-energia-2}.\\
    Now, on the left-hand side of \eqref{eq-energia-2} using the regularity of $u$ and dominated convergence, it  follows that we can pass the limit and get our identity.
     \fin
	\end{dem}\\

To prove that $u$ has an upper bound, we will need the following lemma, which we will prove using Fourier series.
\begin{lem}\label{3}
Let $u$ a $2\pi$-periodic function, then
$$
 \int\limits_{0}^{2\pi} (\triangle)^{s}u(t)u'(t)dt= 0.
$$
\end{lem}

\begin{dem}
Now we write
\begin{equation*}
  u(t)= a_{0} + \sum\limits_{n=1}^{+ \infty} a_{n}\cos(nt) + b_{n}\sin(nt),
\end{equation*}
so,
 \begin{equation*}
  u'(t)=  \sum\limits_{k=1}^{+ \infty} b_{k}k\cos(kt) - a_{k}k\sin(kt),
\end{equation*}
and
\begin{equation*}
  (\triangle)^{s}u(t)=\sum\limits_{n=1}^{+ \infty} a_{n}n^{2s}\cos(nt) + b_{n}n^{2s}\sin(nt),
\end{equation*}

by  orthogonality we have
\begin{equation*}
  <\sin (nt),\sin(kt)> = 0 , n	\neq k
\end{equation*}
\begin{equation*}
  <\cos (nt),\cos(kt)> = 0 , n	\neq k
\end{equation*}
\begin{equation*}
  <\cos (nt),\sin(kt)> = 0.
\end{equation*}
Hence,
\begin{equation*}
 \int\limits_{0}^{2\pi} (\triangle)^{s}u(t)u'(t)dt=\\
\end{equation*}
\begin{equation*}
   \qquad \sum\limits_{n=1}^{+ \infty}\int\limits_{0}^{2\pi} (a_{n}n^{2s}\cos(nt) + b_{n}n^{2s}\sin(nt))( b_{n}n\cos(nt) - a_{n}n\sin(nt))
\end{equation*}
\begin{eqnarray}\label{zero}
   &=&  \sum\limits_{n=1}^{+ \infty}\left(\int\limits_{0}^{2\pi} a_{n}b_{n}n^{2s + 1}\cos^{2}(nt) + \int\limits_{0}^{2\pi}(b_{n}^{2}- a_{n}^{2})n^{2s + 1}\cos(nt)\sin(nt)\right.\nonumber\\
    &&\left. \qquad\qquad\qquad\qquad\qquad\qquad\qquad\qquad\qquad- \int\limits_{0}^{2\pi} a_{n}b_{n}n^{2s + 1}\sin^{2}(nt)\right)\nonumber \\
   &=& \sum\limits_{n=1}^{+ \infty} \int\limits_{0}^{2\pi} a_{n}b_{n}n^{2s + 1}( \cos^{2}(nt) - \sin^{2}(nt))\nonumber\\
   &=& \sum\limits_{n=1}^{+ \infty} \int\limits_{0}^{2\pi} a_{n}b_{n}n^{2s + 1} \cos(2nt) =0.
\end{eqnarray}
\fin
\end{dem}

As mentioned in the introduction, the fact that we can't establish that $u'$ is uniformly bounded, we need to use local regularity theory,
the following Lemma is proved in \cite{felmer}, (see Theorem 3.1).

\begin{lem}\label{regu}
Let $s>\frac{1}{2}$, for $\gamma\in(0,1)$ we have $(\triangle)^{1 -s} : C^{2,\gamma}(\mathbb{R}) \to C^{2s +\gamma}(\mathbb{R})$ is continuous, i.e.
\begin{equation*}
  \|(\triangle)^{1 -s} w\|_{C^{2s + \gamma}(\mathbb{R})} \leq C\|w\|_{C^{2,\gamma}(\mathbb{R})}.
\end{equation*}
\end{lem}
The following local regularity result, we use ideas of \cite{che}, see also \cite{silvestre}, we give it here for completeness.
Moreover, the the basic ideas are the key elements in Propositions 2.8 and 2.9 of \cite{silvestre} also use in this paper.

\begin{lem}\label{local}
 Let $u\in L^{\infty}(\mathbb{R})$ be a solution of
\begin{equation*}
  (\triangle)^{s}u= f \qquad \text{in} \quad (x_{0} - \delta, x_{0} + \delta),
\end{equation*}
with $\delta >0$. Then there exist $ \gamma >0$ and $C^*$ such that $u\in C^{2s + \gamma}_{loc}(\mathbb{R})$. Moreover,
\begin{equation*}
  \|u\|_{C^{2s + \gamma} [x_{0} - \frac{\delta}{2}, x_{0} + \frac{\delta}{2}]}\leq C^*(\|f\|_{C^{\gamma}(x_{0} - \delta, x_{0} + \delta)} + \|u\|_{L^{\infty}(\mathbb{R})}).
\end{equation*}
\end{lem}

\begin{dem}
Let $w$ be a solution of
\begin{equation*}
  w''= \eta f \qquad in\quad \mathbb{R},
\end{equation*}
where $\eta\in C^{\infty}(\mathbb{R})$ such that $\eta\equiv 0$ outside $(x_{0} - \delta, x_{0} + \delta)$ and $\eta\equiv 1$ in $[x_{0} - \frac{3\delta}{4}, x_{0} + \frac{3\delta}{4}]$. So we get
\begin{eqnarray}\label{cw}
  \|w\|_{C^{2,\alpha}(\mathbb{R})} &\leq& C\| \eta f \|_{C^{\alpha}(\mathbb{R})}\nonumber \\
   &\leq& C\|f\|_{C^{\alpha}(x_{0} - \delta, x_{0} + \delta)}.
\end{eqnarray}
Then, since $(\triangle)^{s}((\triangle)^{1 - s}w)= w''$ we have
\begin{equation*}
  (\triangle)^{s}(u - (\triangle)^{1 - s}w)=0 \quad in \quad [x_{0} - 3\delta/4, x_{0} + 3\delta/4],
\end{equation*}
we can use Theorem $1.1$ of \cite{caff}, to obtain that there exist $\gamma$ such that
\begin{eqnarray*}
  \|u - (\triangle)^{1 - s}w\|_{C^{2s+ \gamma}[x_{0} - \frac{3\delta}{4}, x_{0} + \frac{3\delta}{4}]}&\leq& C\|u - (\triangle)^{1 - s}w\|_{L^{\infty}(\mathbb{R})}\\
                                       &\leq& C (\|u\|_{L^{\infty}(\mathbb{R})} + \|(\triangle)^{1 - s}w\|_{L^{\infty}(\mathbb{R})}).
\end{eqnarray*}
Hence, using the Lemma \ref{regu} and (\ref{cw}) we get
\begin{eqnarray*}
    \|u\|_{C^{2s+ \gamma}[x_{0} - \frac{3\delta}{4}, x_{0} + \frac{3\delta}{4}]} &\leq & C (\|(\triangle)^{1 - s}w\|_{C^{2s+\gamma}(\mathbb{R})} + \|u\|_{L^{\infty}(\mathbb{R})}) \\
                                                                             &\leq & C(\|f\|_{C^{\gamma}(x_{0} - \delta, x_{0} + \delta)} + \|u\|_{L^{\infty}(\mathbb{R})}).
\end{eqnarray*}
\fin
\end{dem}

Now use again interpolation inequality to get a local regularity result with drift term.

\begin{teo}\label{tloc}
Let $u$ be a solution of
\begin{equation*}
  (\triangle)^{s}u + cu'= f \qquad \text{in} \quad (x_{0} - \delta, x_{0} + \delta),
\end{equation*}
with $\delta >0$. Then there exists $ \gamma >0$ and $\bar C>0$ such that
\begin{equation*}
  \|u\|_{C^{2s + \gamma} [x_{0} - \frac{\delta}{2}, x_{0} + \frac{\delta}{2}]}\leq \bar C(\|f\|_{C^{\gamma}(x_{0} - \delta, x_{0} + \delta)} + \|u\|_{L^{\infty}(\mathbb{R})}).
\end{equation*}
\end{teo}
\begin{dem} Now, by the Lemma \ref{local},
\begin{equation}\label{rd}
  \|u\|_{C^{2s + \gamma} [x_{0} - \frac{\delta}{2}, x_{0} + \frac{\delta}{2}]} \leq C^*(\|f\|_{C^{\gamma}(x_{0} - \delta, x_{0} + \delta)} + \|u\|_{C^{1,\gamma}(x_{0} - \delta, x_{0} + \delta)} + \|u\|_{L^{\infty}(\mathbb{R})}).
  \end{equation}

Using now the Interpolation inequalities [\cite{interpol}, Theorems 3.2.1] together with the fact that $2s > 1$, we obtain that, for any $\epsilon > 0$, there is a positive constant $C = C(\gamma, \epsilon, s)$, such that
\begin{equation}\label{rp}
 \|u\|_{C^{1, \gamma}(x_{0} - \delta, x_{0} + \delta)}\leq \epsilon  \|u\|_{C^{2s+\gamma}(x_{0} - \delta, x_{0} + \delta)} + C\|u\|_{C(x_{0} - \delta, x_{0} + \delta)}.
\end{equation}

Hence, choosing $\epsilon= \frac{1}{2C^*}$ we get out results from $(\ref{rp})$ and  $(\ref{rd})$.
\fin
\end{dem}

\begin{lem}\label{harnack}(Interior Harnack inequality).
Let $v$ be a classical solution of
\begin{equation*}
  (\triangle)^{s}v + \delta^{2s-1}cv'= f \qquad\text{in}\quad (-1, 1),
\end{equation*}
and $v\geq0$ in $\mathbb{R}$ with $f\in L^{\infty}(-1,1)\cap C(-1,1)$ and $\delta\in (0,1)$.
Then there exists $C_{0}>0$ independent of $v$ and $\delta$ such that
\begin{equation*}
  \sup\limits_{(-\frac{1}{2},\frac{1}{2})} v \leq C_{0}( \inf\limits_{(-1,1)} v + \|f\|_{L^{\infty}(-1,1)}).
\end{equation*}
\end{lem}

For the proof we quote \cite{DQT} where a much more general equations are consider, including zero order term and drift term.
A parabolic version of Harnack inequality with a drift term can be found in \cite{DL}.

\begin{lem}\label{2}
Assume that $(G1)$ and $(G2)$ holds then there exist constants $R_{1} > R_{0} > 0$ such that for each possible $2\pi$-periodic solution $u$ of (\ref{ecu2l}) there exist $t_{0}$, $t_{1} \in [0, 2\pi]$ such that $u(t_{0}) > R_{0}$ and $u(t_{1}) < R_{1}$.
\end{lem}

\begin{dem}
We shall assume that $u$ is a solution of (\ref{ecu2l}) for some fixed $\lambda\in (0,1)$ then
\begin{equation*}
  \int\limits_{0}^{2\pi} (\triangle) ^{s}u(t)dt + c\int\limits_{0}^{2\pi} u'(t)dt - \lambda \int\limits_{0}^{2\pi} g(u(t))dt = \lambda \int\limits_{0}^{2\pi} e(t)dt.
\end{equation*}
But, for Lemma \ref{perizero} and using the $2\pi$-periodicity of $u'$ gives
\begin{equation}\label{igual}
  \frac{1}{2\pi} \int\limits_{0}^{2\pi} g(u(t))dt + \bar{e}=0.
\end{equation}
We now notice that with assumption $(G2)$, we get that there exists  $R_{0}>0$ such that
$$
g(x)>0 \;\; \text{and} \;\;\; g(x) + \bar{e} >0,
$$
whenever $0<x\leq R_{0}$. Therefore, if $0<u(t)\leq R_{0}$ for all $t\in[0, 2\pi]$, we obtain $ g(u(t)) + \bar{e} >$ for those $t$ and hence
\begin{equation*}
  \frac{1}{2\pi} \int\limits_{0}^{2\pi} g(u(t))dt + \bar{e} > 0
\end{equation*}
a contradiction to (\ref{igual}), thus there exist $t_{0}$ such that $u(t_{0}) > R_{0}$. On the
other hand, assumption $(G1)$  implies the existence of some $R_{1} > R_{0}$ such that
\begin{equation*}
  g(u) + \bar{e} < 0,
\end{equation*}
whenever $u \geq R_{1}$. Then, if $u(t) \geq R_{1}$ for all $t\in[0,2\pi]$, gives  $g(u(t)) + \bar{e} < 0$ for those $t$ and
\begin{equation*}
  \frac{1}{2\pi} \int\limits_{0}^{2\pi} g(u(t))dt + \bar{e} < 0
\end{equation*}
a contradiction to (\ref{igual}), thus there exist $t_{1}$ such that $u(t_{1}) < R_{1}$.
\fin
\end{dem}

\begin{lem}\label{cotasup}
Let $u$ a solution positive $2\pi$-periodic of (\ref{ecu2l}) and assume $(G1)$ then there exist constant $R$ such that
$$
0<u(t) < R.
$$
Moreover we have that there exist $C>0$ such that
\begin{equation}\label{deriv}
 \|u'\|_{L^{2}(0,2\pi)} \leq C\|e\|_{L^{2}(0,2\pi)}.
\end{equation}
\end{lem}

\begin{dem}
By Lemma \ref{2}
\begin{eqnarray}\label{cl1}
  u(t)&\leq& u(t_{1}) + \int\limits_{0}^{2\pi} |u'(t)|dt\nonumber  \\
       &\leq& R_{1} + \sqrt{2\pi} \|u'\|_{L^{2}(0, 2\pi)}.
\end{eqnarray}
Now, multiplying (\ref{ecu2l}) for $u'$ and integrating the equation
\begin{equation}\label{laecua}
  \int\limits_{0}^{2\pi} (\triangle) ^{s}u(t)u'(t) dt + c \int\limits_{0}^{2\pi} |u'(t)|^{2}dt = \lambda \int\limits_{0}^{2\pi} g(u(t))u'(t) + \lambda \int\limits_{0}^{2\pi} e(t)u'(t)dt,
\end{equation}
since
\begin{equation*}
  \lambda \int\limits_{0}^{2\pi} g(u(t))u'(t)dt = \lambda \int\limits_{u(0)}^{u(2\pi)} g(w)dw= 0,
\end{equation*}
together with the Lemma \ref{3} and using H\"older inequality gives
\begin{eqnarray*}
  c\|u'\|_{L^{2}(0,2\pi)}^{2} &=& \lambda \int\limits_{0}^{2\pi} e(t)u'(t) \\
   &\leq& C \|e\|_{L^{2}(0,2\pi)}\|u'\|_{L^{2}(0,2\pi)},
\end{eqnarray*}
so that,
\begin{equation*}
 \|u'\|_{L^{2}(0,2\pi)} \leq C\|e\|_{L^{2}(0,2\pi)}
\end{equation*}
which, together with (\ref{cl1}) implies
\begin{equation*}
  u(t) \leq R_{0} + C\|e\|_{L^{2}(0,2\pi)}= R,
\end{equation*}
therefore the result follows.
\fin
\end{dem}

\begin{lem}\label{lemaimport}
Assume that $(G4)$ holds then for each $u$ a solution of $(\ref{ecuag2})$ there exists $C>0$ (independent of $u$) such that
\begin{equation}\label{scota}
  \|u\|_{C^{s + \epsilon/2}[0,2\pi]} \leq C\left(\left(\int\limits_{-2\pi}^{4\pi} |g(u(x))|^{1/(s-\epsilon/2)}dx\right)^{s-\epsilon/2}+ \|u'\|_{L^{2}[-2\pi,4\pi]} + \|e\|_{L^{2}[-2\pi,4\pi]}\right).
\end{equation}
\end{lem}
\begin{dem}
First we find a $L^{1}$ bound for $g$. Notice that $(G4)$ implies that
\begin{equation}\label{gcota}
  |g(u(t))|\leq g(u(t)) + 2au(t) + 2b,
\end{equation}
integrating over $[0,2\pi]$
\begin{equation*}
  \int\limits_{0}^{2\pi} |g(u(t))|dt \leq \int\limits_{0}^{2\pi} g(u(t))dt +2a\int\limits_{0}^{2\pi} u(t)dt + 4b\pi,
\end{equation*}
by Lemma \ref{cotasup} and $(\ref{igual})$ we have
\begin{equation}\label{gl1}
  \int\limits_{0}^{2\pi} |g(u(t))|dt \leq 4\pi(aR + b).
\end{equation}
 Secondly, we now consider the following problem
\begin{equation}\label{aux}
   (\Delta)^s v(t) =g(v(t)) + h(t)\qquad \text{in} \quad \mathbb{R},
\end{equation}
where $h\in L^{2}$. Let $w\in C^{2}(\mathbb{R})$ such that
\begin{equation}\label{w}
  w'(t)=   \int\limits_{-2\pi}^{t} g(v(x)) + h(x)dx, \qquad t\in[-2\pi,4\pi].
\end{equation}
Then, $w$ is a solution of $$w''(t)= g(v(t)) + h(t),$$ and since $(\triangle)^{s}((\triangle)^{1 - s}w)= w''$ we have
\begin{equation*}
  (\triangle)^{s}(v - (\triangle)^{1 - s}w)=0 \quad in \quad [0,2\pi].
\end{equation*}
Using \cite[Theorem 2.7]{CS} we have
\begin{eqnarray*}
  \|v - (\triangle)^{1 - s}w\|_{C^{s + \epsilon/2}(0,2\pi)}&\leq& C\|v - (\triangle)^{1 - s}w\|_{L^{\infty}(0,2\pi)}\\
                                       &\leq& C (\|v\|_{L^{\infty}(\mathbb{R})} + \|(\triangle)^{1 - s}w\|_{L^{\infty}(0,2\pi)}),
\end{eqnarray*}
and as we know that for $\mu=1 - s + \epsilon/2$  (see \cite[Theorem 3.1]{felmer})
\begin{equation}\label{1s}
 \| (\triangle)^{1-s}w\|_{C^{s + \epsilon/2}[0,2\pi]} \leq C \|w\|_{C^{1,\mu}(-2\pi,4\pi)},
\end{equation}
we obtain
\begin{eqnarray}\label{se1}
    \|v\|_{C^{s+\epsilon/2}[0,2\pi]} &\leq & C (\|(\triangle)^{1 - s}w\|_{C^{s+ \epsilon/2}[0,2\pi]} + \|v\|_{L^{\infty}(\mathbb{R})}) \nonumber\\
                                                                             &\leq & C(\|w\|_{C^{1,\mu}(-2\pi,4\pi)} + \|v\|_{L^{\infty}(\mathbb{R})}).
\end{eqnarray}
Now, we want to bound $\|w\|_{C^{1,\mu}(-2\pi,4\pi)}$, it follows from $(\ref{w})$ and $(\ref{gl1})$ we need to estimate $[w']_{C^{\mu}}$. Let $x,y\in(-2\pi,4\pi)$ and using  H\"older inequality we get
\begin{eqnarray}\label{wscota}
  |w'(x) -w'(y)| &=& \left|\int\limits_{x}^{y} g(v(x)) + h(x)dx\right| \nonumber\\
                 &\leq& \left(\left(\int\limits_{x}^{y}|g(v(x))|^{1/(s-\epsilon/2)}dx\right)^{s-\epsilon/2}+ \left(\int\limits_{x}^{y}|h(x)|^{1/(s-\epsilon/2)}dx\right)^{s-\epsilon/2}\right)|x-y|^{1-s+\epsilon/2}\nonumber\\
                 &\leq& \left(\left(\int\limits_{-2\pi}^{4\pi} |g(v(x))|^{1/(s-\epsilon/2)}dx\right)^{s-\epsilon/2} + \left(\int\limits_{-2\pi}^{4\pi} |h(x)|^{1/(s-\epsilon/2)}dx\right)^{s-\epsilon/2}\right)|x-y|^{1-s+\epsilon/2},\nonumber\\
                 &&
\end{eqnarray}
so that $(\ref{se1})$ implies,
\begin{equation}\label{cv}
    \|v\|_{C^{s + \epsilon/2}[0,2\pi]} \leq C\left(\left(\int\limits_{-2\pi}^{4\pi} |g(v(x))|^{1/(s-\epsilon/2)}dx\right)^{s-\epsilon/2} + \|h\|_{L^{2}[-2\pi,4\pi]}\right).
\end{equation}
Finally, let $u$ a solution of $(\ref{ecuag2})$ by the above with  $h(t)= cu'(t) + e(t)$  we have
\begin{equation*}
  \|u\|_{C^{s + \epsilon/2}[0,2\pi]} \leq C\left(\left(\int\limits_{-2\pi}^{4\pi} |g(u(x))|^{1/(s-\epsilon/2)}dx\right)^{s-\epsilon/2}+ \|u'\|_{L^{2}[-2\pi,4\pi]} + \|e\|_{L^{2}[-2\pi,4\pi]}\right).
\end{equation*}
\fin
\end{dem}

\begin{lem}\label{cotainf}
Assume $(G1)$,$(G2)$,$(G3)$ and $(G4)$ then there exists $r \in(0,R_{0})$ such that each $2\pi$-periodic solution of (\ref{ecu2l}) satisfies $u(t) > r$ for all $t \in [0, 2\pi]$.
\end{lem}

\begin{dem}
Let $t^n_{3}$, $t^n_{4}$ be the minimum point and the maximum point of $u_n$ in $[0,2\pi]$, notice that by Lemma \ref{2}  be have $u_n(t^n_4)>R_0$.

Now we assume by contradiction that $u_n(t^n_3) \to 0$ as $n\to \infty$.  First we  suppose that $t^n_{4}<t^n_{3}$,

and multiplying (\ref{ecu2l}) by $u'_n$, then
\begin{equation*}
   \int\limits_{t^n_{4}}^{t^n_{3}} (\triangle) ^{s}u_n(t)u_n'(t) dt + c \int\limits_{t^n_{4}}^{t^n_{3}} |u_n'(t)|^{2}dt - \lambda\int\limits_{t^n_{4}}^{t^n_{3}} g(u_n(t))u_n'(t) = \lambda \int\limits_{t^n_4}^{t^n_{3}} e(t)u_n'(t)dt,
\end{equation*}
so that,
\begin{equation}\label{se}
   \lambda\int\limits_{u_n(t^n_{3})}^{u_n(t^n_{4})} g(w)dw= \int\limits_{t^n_{4}}^{t^n_{3}} (-\triangle) ^{s}u_n(t)u_n'(t) dt - c \int\limits_{t^n_{4}}^{t^n_{3}} |u_n'(t)|^{2}dt + \lambda\int\limits_{t^n_{4}}^{t^n_{3}} e(t)u_n'(t)dt.
\end{equation}
We need to bound the first term in the right hand side, since by (\ref{deriv}) the to order terms in the right hand are bounded,
to get a contradiction with $(G3)$ by the fact that $u_n(t^n_4)>R_0$ and $u_n(t^n_3) \to 0$.

Let as now define $$\delta_n=\sup\{\delta>0\,\,| u_n (t)>\frac{R_0}{4C_0} \mbox{ for  } t \in (t_{4}^{n}-\delta, t_{4}^{n}+\delta)\}.$$
Notice that for $n$ large $u_n(t_{3}^{n})<\frac{R_0}{4C_0}$ therefore $\delta_n$ is finite and
$$\inf_{(t_{4}^{n}-\delta_n, t_{4}^{n}+\delta_n)}u_n=\frac{R_0}{4C_0}.$$

Now we claim that there exists $\delta_{0}>0$ such that $\delta_n>\delta_{0}$.
Suppose the contrary, so there exists a sub-sequence (still denote by $n$) such that $\delta_{n} \to 0$. Define $v_{n}(t)=u_{n}(\delta_{n}t + t_{4}^{n})$, we have that
\begin{equation*}
  (\triangle)^{s}v_{n}(t) + \delta^{2s-1}_{n}cv_n'(t)= \delta^{2s}_{n}h \qquad in \quad (-1,1),
\end{equation*}
where $h:=g(u_{n}(\delta_{n}t + t_{4}^{n})) + e(\delta_{n}t + t_{4}^{n})\in L^\infty(-1,1)$ by the definition of $\delta_n$ .
By the Lemma \ref{harnack},
\begin{equation*}
   \sup\limits_{(-\frac{1}{2},\frac{1}{2})} v_{n} \leq C_{0}( \inf\limits_{(-1,1)} v_{n} + \delta^{2s}_{n}\|h\|_{L^{\infty}(-1,1)}),
\end{equation*}
as $\inf\limits_{(-1,1)} v_{n} = \dfrac{R_{0}}{4C_{0}}$,  we have
\begin{equation*}
  R_{0} <  \sup\limits_{(-\frac{1}{2},\frac{1}{2})} v_{n} \leq C_{0}(\frac{R_{0}}{4C_{0}} + \delta^{2s}_{n}\|h\|_{L^{\infty}(-1,1)}).
\end{equation*}
Taking $\delta_{n}$ small, we obtain a contradiction and the claim follows.

Now we can use Theorem \ref{tloc} and with $\delta=\delta_0 >0$ independent of $u_n$ and get there exist $\bar C>0$ (independent of $n$) such that
\begin{equation*}
   \|u_n\|_{C^{2s + \beta} [t_{4}^{n} - \frac{\delta}{2}, t_{4}^{n} + \frac{\delta}{2}]}\leq \bar C,
\end{equation*}
with $\delta > 0$, this implies that
\begin{equation}\label{cotauprima}
|u_n'(x)| \leq \bar C\quad\mbox{ in }\quad [t_{4}^{n} - \delta/2, t_{4}^{n} + \delta/2].
\end{equation}

Now we simplify the notation and drop the $n$ index to estimate  $\int\limits_{t^n_{4}}^{t^n_{3}} u_n'(x) (-\Delta)^s u_n(x) dx$.

By the Lemma \ref{lemint} we have
\begin{eqnarray}\label{ener}
	  \int\limits_{t_{4}}^{t_{3}} u'(x) (-\Delta)^s u(x) dx &\leq&\frac{c(1,s)}{2} \left( \int_{-\infty}^{+\infty} \frac{(u(t_{3})-u(y))^2}{|t_{3}-y|^{1+2s}} dy    \right.\nonumber\\
     &&  \left. + (1+2s)\int\limits_{t_{4}}^{t_{3}} \int\limits_{-\infty}^{t_{4}} \frac{(u(x)-u(y))^2}{|x-y|^{2+2s}} dy dx\right).
	\end{eqnarray}
Let $\rho>0$ we have
\begin{equation*}
   \int_{-\infty}^{+\infty} \frac{(u(t_{3})-u(y))^2}{|t_{3}-y|^{1+2s}} dy =  \int_{|t_{3} - y|< \rho} \frac{(u(t_{3})-u(y))^2}{|t_{3}-y|^{1+2s}} dy +  \int_{|t_{3} - y|\geq \rho} \frac{(u(t_{3})-u(y))^2}{|t_{3}-y|^{1+2s}} dy,
\end{equation*}
by Lemma \ref{lemaimport} together with $(G2)$,$(G4)$ and $(\ref{gcota})$ we obtain
\begin{eqnarray}\label{c1}
  \int\limits_{|t_{3} - y|< \rho} \frac{(u(t_{3})-u(y))^2}{|t_{3}-y|^{1+2s}} dy  &\leq& C\left(\left(\int\limits_{-2\pi}^{4\pi} |g(u(x))|^{1/(s-\epsilon/2)}dx\right)^{2s-\epsilon} + C^{*}\right)\int_{|t_{3} - y|< \rho} \frac{|t_{3}-y|^{2s+\epsilon}}{|t_{3}-y|^{1+2s}} dy  \nonumber \\
  &\leq& C \left(\left(\int\limits_{-2\pi}^{4\pi} |g(u(x))|^{1 + (1/(s-\epsilon2) -1)}dx\right)^{2s-\epsilon} + C^{*}\right)\int_{|t_{3} - y|< \rho} |t_{3}-y|^{\epsilon -1} dy\nonumber\\
  &\leq& Cg(u(t_{3}))^{(1/(s-\epsilon/2) -1)(2s-\epsilon)} \left(\left(\int\limits_{-2\pi}^{4\pi} g(u(x))dx\right)^{2s-\epsilon} + + C^{*}\right)\nonumber\\
  &\leq& C(g(u(t_{3}))^{(1/(s-\epsilon/2) -1)(2s-\epsilon)} + C^{*})=C(g(u(t_{3}))^{2-2s +\epsilon} + C^{*}),
\end{eqnarray}
where $C^{*}= 4a\pi R + 4b\pi + \|u'\|_{L^{2}[-2\pi,4\pi]} + \|e\|_{L^{2}[-2\pi,4\pi]}$, and by $(\ref{deriv})$ we have $C^{*}< + \infty$. Now, the second integral we know is bounded
\begin{eqnarray}\label{c4}
  \int_{|t_{3} - y|\geq \rho}  \frac{(u(y) -u(t_{3}))^{2}}{|t_{3}-y|^{1+2s}} dy &\leq& 4R^{2} \int_{|t_{3} - y|\geq \rho} \frac{dy}{|t_{3}-y|^{1+2s}}\nonumber\\
   &\leq& C.
\end{eqnarray}

Now from \eqref{cotauprima} we have
\begin{equation*}
  u(x) - u(y) = \int\limits_{y}^{x} u'(\xi)d\xi\leq \bar C (x-y),
\end{equation*}
so that
\begin{eqnarray}\label{c2}
   \int\limits_{t_{4}}^{t_{3}} \int\limits_{t_{4}- \frac{\delta}{2}}^{t_{4}} \frac{(u(x)-u(y))^2}{|x-y|^{2+2s}} dy dx &\leq& \bar C^{2}  \int\limits_{t_{4}}^{t_{3}} \int\limits_{t_{4}- \frac{\delta}{2}}^{t_{4}} \frac{(x -y)^2}{|x-y|^{2+2s}} dy dx \nonumber \\
   &=& \dfrac{\bar C^{2}}{2s -1} \int\limits_{t_{4}}^{t_{3}} -(x-t_{4})^{-2s +1} + \left(x-t_{4}+ \frac{\delta}{2}\right)^{-2s +1} dx\nonumber \\
   &=& \dfrac{\bar C^{2}}{(2s -1)(2 -2s)} \left( -(t_{3}-t_{4})^{-2s +2} + \left(t_{3}-t_{4} + \frac{\delta}{2}\right)^{-2s +2}\right. \nonumber\\
    && \qquad\qquad\qquad\qquad\qquad\qquad\qquad\qquad\qquad \left.- \frac{\delta}{2}^{-2s +2}\right)\nonumber\\
   &\leq& \dfrac{\bar C^{2}}{(2s -1)(2 -2s)} \left(t_{3}-t_{4} + \frac{\delta}{2}\right)^{-2s +2}.
\end{eqnarray}

To analyze the same integral when $y < t_{4}- \delta$ and $t_{4} < x< t_{3}$ we have
\begin{eqnarray}\label{c3}
  \int\limits_{t_{4}}^{t_{3}} \int\limits_{- \infty}^{t_{4} - \frac{\delta}{2}} \frac{(u(x)-u(y))^2}{|x-y|^{2+2s}} dy dx  &\leq&   4\|u\|^{2}_{L_{\infty}}\int\limits_{t_{4}}^{t_{3}} \int\limits_{- \infty}^{t_{4} - \frac{\delta}{2}} \frac{dy}{|x-y|^{2+2s}} dx \nonumber\\
   &\leq& \frac{4R^{2}}{1+2s}\int\limits_{t_{4}}^{t_{3}}- \left(x - t_{4} + \frac{\delta}{2}\right)^{-1 -2s}dx \nonumber\\
   &=& \frac{4R^{2}}{(1 +2s)(2s)}\left(\left(t_{3}-t_{4} + \frac{\delta}{2}\right)^{-2s} - \frac{\delta}{2}^{-2s}\right)\nonumber\\
   &\leq& \frac{4R^{2}}{(1 +2s)(2s)}\left(t_{3}-t_{4} + \frac{\delta}{2}\right)^{-2s}.
\end{eqnarray}
Hence, \eqref{ener} $(\ref{c1})$, $(\ref{c4})$, $(\ref{c2})$ and $(\ref{c3})$  gives that there exist $M_0>0$ independent of $n$ such that
\begin{equation*}
   g(u_n(t^{n}_{3}))^{2s-2-\epsilon}\int\limits_{t^n_{4}}^{t^n_{3}} u_n'(x) (-\Delta)^s u_n(x) dx \leq  M_{0}.
\end{equation*}
Thus, as mentioned above, this inequality implies a contradiction and the result follows.

\fin
\end{dem}\\
We now prove the following existence result for the $2\pi$-periodic solutions of equation (\ref{ecuag2}), we shall use ideas the Continuation Theorem of \cite{mawhin, gain}.\\

\begin{deu}

By the definition of $R_0$ and $R_1$ in Lemma \ref{2} we have
\begin{equation}\label{ct}
 g(a) >0\;\; \text{and} \; g(a) + \bar{e}> 0\;\;  \text{if} \quad 0< a\leq r < R_{0}
\end{equation}
and
\begin{equation}\label{ctt}
  g(a) + \bar{e}< 0\;\;  \text{if} \quad a\geq R> R_{1}.
\end{equation}

Using Proposition $\ref{l1}$ we can define a map $K:C^{\alpha}_{2\pi}(\mathbb{R}) \to C^{\alpha}_{2\pi}(\mathbb{R})$ by $K(z)= u$ where $u$ is a solution of
\begin{equation*}
  (\triangle)^{s}u(t) + cu'(t)= z(t),
\end{equation*}
where $K$ is compact.\\

Now, let us define the map $N: C^{\alpha}_{2\pi}(\mathbb{R}) \to C^{\alpha}_{2\pi}(\mathbb{R})$ by
\begin{equation}\label{defnn}
  Nu= g(u(\cdot)) + e(\cdot).
\end{equation}

 Define the continuous projectors $Q:  C^{\alpha}_{2\pi}(\mathbb{R}) \to C^{\alpha}_{2\pi}(\mathbb{R})$ by the constant function
 \begin{equation*}
   Qy= \frac{1}{2\pi}   \int\limits_{0}^{2\pi} y(t)dt.
 \end{equation*}

Let $\Omega= \{ u \in C^{2s +\alpha}_{2\pi}(\mathbb{R}): r <u(t) <R,\;  t\in [0,2\pi] \}$  define one parameter family of problems
\begin{equation*}
   u=K((1-\lambda)QNu + \lambda Nu)  \;\; \lambda \in [0,1].
\end{equation*}
Explicitly,
\begin{equation}\label{homo}
   (\triangle)^{s}u(t) + cu'(t)= (1-\lambda)QNu + \lambda Nu.
\end{equation}
For $\lambda \in [0, 1]$, observe that we have by Lemma \ref{perizero}
\begin{equation*}
  \frac{1}{2\pi}\int\limits_{0}^{2\pi} g(u(t)) + e(t) dt=0.
\end{equation*}
Therefore for all $\lambda\in (0, 1]$, problem $(\ref{ecu2l})$ and problem $(\ref{homo})$ are equivalent. Hence, Lemmas \ref{cotasup} and \ref{cotainf} implies  $(\ref{homo})$ does not have a solution on $\partial\Omega \times (0,1]$. For, $\lambda=0$ $(\ref{homo})$ is equivalent to the problem
\begin{equation}\label{qn}
   (\triangle)^{s}u(t) + cu'(t)=  \frac{1}{2\pi}\int\limits_{0}^{2\pi} g(u(t)) + e(t) dt,
\end{equation}
 then, applying $Q$ to both members of this equation and by Lemma \ref{perizero}, we obtain
 \begin{eqnarray*}
   QNu=0, &&   (\triangle)^{s}u(t) + cu'(t)= 0.
 \end{eqnarray*}
 Multiplying by $u'$ the second of those equations, using the Lemma \ref{3} and integrating we have
\begin{equation*}
  \int\limits_{0}^{2\pi} (\triangle)^{s}u(t) u'(t) dt + c\int\limits_{0}^{2\pi} (u'(t))^{2}dt = 0,
\end{equation*}
so that
\begin{equation*}
  \| u' \|_{L^{2}(0, 2\pi)} = 0,
\end{equation*}
this implies $u'=0$ hence $u$ is constant, and we know that the constant solutions of $QNa =0$ satisfy the inequality $r<a<R$. Thus we have proof that
$(\ref{homo})$ has no solution on $\partial\Omega\times [0,1]$. Therefore, the $deg(I- K((1-\lambda)QN + \lambda N), \Omega, 0)$  is well defined  for all $\lambda\in[0,1]$ and by homotopy invariant of the degree we have
 \begin{equation*}
   deg(I - KN, \Omega, 0)= deg(I - KQN, \Omega, 0)=deg(I - KQN, \Omega\cap E, 0),
 \end{equation*}
where $E\subset C^{2s +\alpha}_{2\pi}(\mathbb{R}):$ is the one dimensional space of constants maps, and the last equality is since all solution of \eqref{qn} are constants maps as proved above.

From here, we use \eqref{ct} and \eqref{ctt} and basic degree properties to get $deg(I - KQN, \Omega\cap E, 0)\neq 0$
Thus, we can conclude equation $(\ref{ecuag2})$ has at least one $2\pi$-periodic classical solution.
\fin
\end{deu}

\section{Bifurcation from infinity and multiplicity of the solutions}
In this section, we discuss a multiplicity result.  We find the existence of a continuum of positive solutions, bifurcating from infinity this together with our previos results
will give a multiplicity of solutions. This result is based on ideas from \cite{mawhin} and \cite{schmitt}.  Here we will use the notation of the previous section.\\

The eigenvalue problem
\begin{equation*}
  u= \mu K u
\end{equation*}
has associated  eigenvalue $\mu=0$ the constant eigenfunction $u\equiv 1$. Conversely, periodic eigenfunctions associated with $\mu=0$ are necessarily constant, (see proof the Theorem 1.4), therefore  $\mu=0$ is a simple eigenvalue.\\
We  want to find positive $2\pi$-periodic solutions of the equation
\begin{equation}\label{bifur}
   (\triangle)^{s}u(x) + cu'(x) + \mu u = G(u) + e(t).
\end{equation}
We will assume that continuous functions $G: (0, + \infty) \to [0, \infty)$ and $e\in  C^{\alpha}_{2\pi}(\mathbb{R})$ satisfy the following  conditions
\begin{enumerate}[(H1)]
  \item $\lim\limits_{t  \to 0^{+}} G(t) = + \infty$,
  \item $\lim\limits_{t  \to + \infty} G(t) = 0$,
  \item $ \bar{e} >0 $,
  \item $\int\limits_{0}^{1} G(t) dt= + \infty$.
\end{enumerate}
We that we have the following result.

\begin{teo}\label{bifu}
Assume that conditions $(H1)$,$(H2)$,$(H3)$, and $(H4)$ are satisfied. Then there exists $\eta > 0$ such that the following holds:
\begin{itemize}
  \item equation (\ref{bifur}) hast at least one positive solution $u$ for $0\leq \mu$.
  \item equation (\ref{bifur}) hast at least two positive  $2\pi$-periodic solutions $u$ for $- \eta \leq \mu < 0$.
\end{itemize}
\end{teo}
\begin{dem}
We now take $g(u)= G(u) - \mu u$, so (\ref{bifur}) is of the form (\ref{ecuag2}), and satisfying the following conditions:
\begin{enumerate}[(H1')]
  \item $\lim\limits_{t  \to 0^{+}} g(t) = + \infty$
  \item $\limsup\limits_{t  \to + \infty} [g(t) + \bar{e}]< 0$  for $\mu \geq 0$
  \item  $\int\limits_{0}^{1} g(t) dt= + \infty$.
\end{enumerate}
Therefore the results of Theorem \ref{teoulti} are valid for equation (\ref{bifur}) when $\mu\geq 0$. Then, by the continuity of degree define in Theorem
\ref{teoulti} there exists $\eta>0$ such that for $-\eta \leq \mu < 0$ that degree is not trivial. So, there exists $u$ a solution for (\ref{bifur}) with $-\eta \leq \mu < 0$.
Now, by $(H2)$ $N(u)=o(\|u\|)$ at $u= + \infty$ then, the fundamental Theorem on bifurcation from infinity from a simple eigenvalue implies the existence of a continuum $\mathcal{C}_{\infty}$ of positive solution $(\mu, u)$ bifurcating from infinity at $\mu =0$, since the solutions for $\mu \geq 0$ are bounded, the bifurcation is for the left side.
\fin
\end{dem}\\

{\bf Acknowledgements}
 A. Q. was partially supported by FONDECYT Grant \# 1190282 and Programa Basal, CMM. U. de Chile.
L. C. was partially supported by ANID \# 21191475 and by "PIIC de la Dirección de Postgrado y Programas de la UTFSM".

\end{document}